\documentclass{amsart}
\usepackage{pstricks,pst-node,amssymb}

\title{Chains in the Bruhat order}

\author{Alexander Postnikov and Richard P.~Stanley}

\address{Department of Mathematics, M.I.T., Cambridge, MA 02139}
\email{apost@math.mit.edu}
\email{rstan@math.mit.edu}


\thanks{A.P.\ was supported in part by National Science Foundation 
grant DMS-0201494 and by Alfred P.\ Sloan Foundation research fellowship}
\thanks{R.S.\ was supported in part by National Science Foundation 
grant DMS-9988459}

\date{September~7, 2004; minor updates on February~07, 2005}

\numberwithin{equation}{section}

\theoremstyle{plain}
\newtheorem{theorem}{Theorem}[section]
\newtheorem{proposition}[theorem]{Proposition}
\newtheorem{lemma}[theorem]{Lemma}
\newtheorem{corollary}[theorem]{Corollary}
\newtheorem{exercise}[theorem]{Exercise}

\theoremstyle{definition}

\newtheorem{conjecture}[theorem]{Conjecture}
\newtheorem{example}[theorem]{Example}

\theoremstyle{remark}
\newtheorem{remark}[theorem]{Remark}

\def\R{\mathbb{R}}
\def\Z{\mathbb{Z}}
\def\N{\mathbb{N}}
\def\C{\mathbb{C}}
\def\Q{\mathbb{Q}}
\def\T{\mathcal{T}}
\def\O{\mathcal{O}}

\def\L{\mathcal{L}}

\def\h{\mathfrak{h}}

\def\wt{\mathrm{wt}}

\def\Vol{\mathrm{Vol}}
\def\id{\mathit{id}}
\def\S{\mathfrak{S}}

\def\D{\mathfrak{D}}
\def\CC{\mathcal{C}}
\def\code{\mathrm{code}}

\def\I{\mathcal{I}}
\def\J{\mathcal{J}}

\def\ddy{{\partial}/{\partial y}}
\def\ddx{{\partial}/{\partial x}}

\def\Sym{\mathit{Sym}}
\def\<{\left<}
\def\>{\right>}
\def\({\left(}
\def\){\right)}
\def\CT{\mathrm{CT}}
\def\ch{\mathit{ch}}
\def\Inv{\mathrm{Inv}}
\def\P{\mathcal{P}}

\def\H{\mathcal{H}}

\begin{document}

\begin{abstract}
We study a family of polynomials whose values express degrees of
Schubert varieties in the generalized complex flag manifold $G/B$.
The polynomials are given by weighted sums over saturated chains
in the Bruhat order.   We derive several explicit formulas for these 
polynomials, and investigate their relations with Schubert polynomials, 
harmonic polynomials, Demazure characters, and
generalized Littlewood-Richardson coefficients.
In the second half of the paper, we concern with the case of
to the classical flag manifold of Lie type $A$ and discuss 
related combinatorial objects:
flagged Schur polynomials, $312$-avoiding permutations, generalized 
Gelfand-Tsetlin polytopes, the inverse Schubert-Kostka matrix,
parking functions, and binary trees.
\end{abstract}

\maketitle

\section{Introduction}
\label{sec:intro}

The complex generalized flag manifold $G/B$ embeds into projective space
$\mathbb{P}(V_\lambda)$, for an irreducible representation $V_\lambda$ of
$G$.  The degree of a Schubert variety $X_w\subset G/B$ in this embedding is a
polynomial function of $\lambda$.  The aim of this paper is to study
the family of polynomials $\D_w$ in $r=\mathrm{rank}(G)$ variables that express 
degrees of Schubert varieties.  According to
Chevalley's formula~\cite{Chev}, also known as Monk's rule in type $A$, these polynomials
are given by weighted sums over saturated chains from $id$ to $w$ 
in the Bruhat order on the Weyl group.
These weighted sums over saturated chains appeared in
Bernstein-Gelfand-Gelfand~\cite{BGG} and
in Lascoux-Sch\"utzenberger~\cite{LS2}.
Stembridge~\cite{Ste} recently investigated these sums in the case when
$w=w_\circ$ is the longest element in the Weyl group.
The value $\D_w(\lambda)$ is also equal to the leading coefficient in
the dimension of the Demazure modules $V_{k\lambda,w}$,
as $k\to\infty$.

The polynomials $\D_w$ are dual to the Schubert polynomials $\S_w$ 
with respect a certain natural pairing on the polynomial ring.  
They form a basis in the space of $W$-harmonic polynomials.  
We show that Bernstein-Gelfand-Gelfand's results~\cite{BGG} easily imply two
different formulas for the polynomials $\D_w$.
The first ``top-to-bottom'' formula starts with the top 
polynomial $\D_{w_\circ}$, which is given by the Vandermonde 
product. 
The remaining polynomials $\D_w$ are obtained from $\D_{w_\circ}$
by applying differential operators associated with Schubert polynomials.
The second ``bottom-to-top'' formula starts with $\D_{\id} = 1$.
The remaining polynomials $\D_w$ are obtained from $\D_{\id}$ by
applying certain integration operators.  Duan's recent result~\cite{Du} 
about degrees of Schubert varieties can be deduced from the bottom-to-top
formula.

Let $c_{u,v}^w$ be the generalized Littlewood-Richardson coefficient
defined as the structure constant of the cohomology ring of $G/B$
in the basis of Schubert classes.  The coefficients $c_{u,v}^w$
are related to the polynomials $\D_w$ in two different ways.
Define a more general collection of polynomials $\D_{u,w}$
as sums over saturated chains from $u$ to $w$ in the Bruhat order 
with similar weights.  (In particular, $\D_w = \D_{\id,w}$.)  
The polynomials $\D_{u,w}$ extend the $\D_w$ in the same way
as the skew Schur polynomials extend the usual Schur polynomials.
The expansion coefficients of $\D_{u,w}$ in the basis of $\D_v$'s
are exactly the generalized Littlewood-Richardson coefficients:
$\D_{u,w} = \sum_v c_{u,v}^w\,\D_v$.
On the other hand, we have $\D_w(y+z) = 
\sum_{u,v} c_{u,v}^w\, \D_u(y)\,\D_w(z)$,
where $\D_w(y+z)$ denote the polynomial in 
pairwise sums of two sets $y$ and $z$ of variables.

We pay closer attention to the Lie type $A$ case. 
In this case, the Weyl group is the symmetric group $W=S_n$.
Schubert polynomials for vexillary permutations, i.e., $2143$-avoiding 
permutations, are known to be given by flagged Schur polynomials.
From this we derive a more explicit formula for the polynomials
$\D_w$ for $3412$-avoiding permutations $w$ and, in particular,
an especially nice determinant expression for $\D_w$ in the case when 
$w$ is $312$-avoiding.

It is well-known that the number of $312$-avoiding permutations in $S_n$
is equal to the Catalan number
$C_n=\frac 1{n+1}\binom{2n}{n}$.  Actually, these permutations are exactly
the Kempf elements studied by Lakshmibai~\cite{Lak}  
(though her definition is quite different).
We show that the characters $ch(V_{\lambda,w})$  of 
Demazure modules for $312$-avoiding permutations
are given by flagged Schur polynomials.  (Here flagged Schur polynomials
appear in a different way than in the previous paragraph.)
This expression can be geometrically interpreted in terms of
generalized Gelfand-Tsetlin polytopes $\P_{\lambda,w}$ studied
by Kogan~\cite{Kog}.
The Demazure character $ch(V_{\lambda,w})$ equals a certain sum over lattice
points in $\P_{\lambda,w}$, and thus, the value $\D_w(\lambda)$ equals 
the normalized volume of  $\P_{\lambda,w}$. 
The generalized Gelfand-Tsetlin polytopes $\P_{\lambda,w}$ are related 
to the toric degeneration of Schubert varieties $X_w$ constructed 
by Conciulea and Lakshmibai~\cite{GL}.

One can expand Schubert polynomials as nonnegative sums of monomials using
RC-graphs.  We call the matrix $K$ of coefficients in these expressions the
{\it Schubert-Kostka matrix,} because it extends the usual Kostka matrix.
It is an open problem to find a subtraction-free
expression for entries of the inverse Schubert-Kostka matrix $K^{-1}$.
The entries of $K^{-1}$ are exactly the coefficients of monomials in the 
polynomials $\D_w$ normalized by a product of factorials.  
On the other hand, the entries of $K^{-1}$ are also the expansion coefficients 
of Schubert polynomials in terms of standard 
elementary monomials.  We give a simple expression for entries of $K^{-1}$ 
corresponding to 312-avoiding permutations and 231-avoiding permutations.  
Actually, these special entries are  always equal to $\pm 1$, or $0$.

We illustrate our results by calculating the polynomial $\D_w$
for the long cycle $w=(1,2,\dots,n)\in S_n$ in five different ways.
First, we show that $\D_w$ equals a sum over parking functions.
This polynomial appeared in Pitman-Stanley~\cite{SP} as the volume of a certain polytope.  
Indeed, the generalized Gelfand-Tsetlin polytope $\P_{\lambda,w}$ 
for the long cycle $w$, which is a 312-avoiding permutation, 
is exactly the polytope studied in~\cite{SP}.
Then the determinant formula leads to another
simple expression for $\D_w$ given by a sum of $2^n$ monomials.
Finally, we calculate $\D_w$ by counting saturated chains in the Bruhat order
and obtain an expression for this polynomial as a sum over binary trees.

The general outline of the paper follows.
In Section~\ref{sec:notations}, we give basic notation related
to root systems. In Section~\ref{sec:schubert-calc}, we recall
classical results about Schubert calculus for $G/B$.
In Section~\ref{sec:deg-schub-var}, we define the polynomials
$\D_w$ and $\D_{u,w}$ and discuss their geometric meaning.  
In Section~\ref{sec:D-pairing}, we discuss the 
pairing on the polynomial ring and harmonic polynomials.
In Section~\ref{sec:expr-for-D}, we prove the top-to-bottom and
the bottom-to-top formulas for the polynomial $\D_w$ and give several 
corollaries.  In particular, we show how these polynomials are
related to the generalized Littlewood-Richardson coefficients.
In Section~\ref{sec:examp-duan}, we give several examples and
deduce Duan's formula.
In Section~\ref{sec:K-th-Demazure}, we recall a few facts about 
the K-theory of $G/B$.
In Section~\ref{sec:asymtotic}, we give a simple proof
of the product formula for $\D_{w_\circ}$.
In Section~\ref{sec:permanent}, we mention a formula for the 
permanent of a certain matrix.
The rest of the paper is concerned with the type $A$ case.
In Section~\ref{sec:A-Schubert-poly}, we recall Lascoux-Sch\"utzenberger's
definition of Schubert polynomials.
In Section~\ref{sec:degree-type-A}, we specialize the results of the
first half of the paper to type $A$.
In Section~\ref{sec:flagged-Schur}, we discuss flagged Schur polynomials,
vexillary and dominant permutations,
and give a simple formula for the polynomials $\D_w$, 
for 312-avoiding permutations.
In Section~\ref{sec:demazure-312}, we give a simple proof of the fact
that Demazure characters for 312-avoiding permutations are given 
by flagged Schur polynomials.
In Section~\ref{sec:generalized-GT}, we interpret this
claim in terms of generalized Gelfand-Tsetlin polytopes.
In Section~\ref{sec:Kostka}, we discuss the inverse of the
Schubert-Kostka matrix.
In Section~\ref{sec:parking}, we discuss the special case of the long cycle
related to parking functions and binary trees.

\smallskip
{\sc Acknowledgments:} \
We thank V.~Lakshmibai for helpful discussions and to Arun Ram for
help with references.

\section{Notations}
\label{sec:notations}

Let $G$ be a complex semisimple simply-connected Lie group.
Fix a Borel subgroup $B$ and a maximal torus $T$ such that 
$G\supset B\supset T$.
Let $\h$ be the corresponding Cartan subalgebra
of the Lie algebra $\mathfrak{g}$ of $G$, and  let $r$ be its rank.
Let $\Phi\subset \h^*$ denote the corresponding {\it root system}. 
Let $\Phi^+\subset \Phi$ be the set 
of positive roots corresponding to our choice of $B$. 
Then $\Phi$ is the disjoint union of $\Phi^+$ and $\Phi^- = -\Phi^+$.
Let $V\subset\h^*$ be the linear space over $\Q$ spanned by $\Phi$.
Let $\alpha_1,\dots,\alpha_r\in\Phi^+$ be the associated set of 
{\it simple roots}.  They form a basis of the space $V$.
Let $(x,y)$ denote the scalar product on $V$ induced by
the Killing form.  
For a root $\alpha\in\Phi$, the corresponding {\it coroot\/} is given 
by $\alpha^\vee = 2\alpha/(\alpha,\alpha)$.  The collection of coroots
forms the dual root system $\Phi^\vee$.

The {\it Weyl group\/} $W\subset\mathrm{Aut}(V)$ 
of the Lie group $G$ is generated by the reflections 
$s_{\alpha}: y \mapsto y - (y,\alpha^\vee)\,\alpha$,
for $\alpha\in\Phi$ and $y\in V$.
Actually, the Weyl group $W$ is generated by 
{\it simple reflections\/} $s_1,\dots,s_r$ corresponding 
to the simple roots, $s_i = s_{\alpha_i}$, subject to the 
{\it Coxeter relations:}
$(s_i)^2=1$ and $(s_i s_j)^{m_{ij}}=1$,
where $m_{ij}$ is half the order of the dihedral subgroup 
generated by $s_i$ and $s_j$.

An expression of a Weyl group element $w$ as a product 
of generators $w=s_{i_1}\cdots s_{i_l}$
of minimal possible length $l$ is called a {\it reduced decomposition\/}
for $w$. Its length $l$ is called the {\it length\/} of $w$
and denoted $\ell(w)$.
The Weyl group $W$ contains a unique {\it longest element\/} $w_\circ$
of maximal possible length $\ell(w_\circ)=|\Phi^+|$.

The {\it Bruhat order\/} on the Weyl group $W$ is the partial order 
relation ``$\leq $'' which is the transitive closure of the 
following covering relation: $u\lessdot w$, for $u,w\in W$, whenever
$w=u \, s_{\alpha}$, for some $\alpha\in\Phi^+$, and $\ell(u)=\ell(w)-1$.
The Bruhat order has the unique minimal element $\id$ and the unique
maximal element $w_\circ$.
This order can also be characterized, as follows.
For a reduced decomposition $w=s_{i_1}\cdots s_{i_l}\in W$ and $u\in W$,
$u\leq  w$ if and only if there exists a reduced decomposition
$u=s_{j_1}\cdots s_{j_s}$ such that $j_1,\dots,j_s$ is a subword
of $i_1,\dots,i_l$.

Let $\Lambda$ denote the {\it weight lattice\/} 
$\Lambda=\{\lambda\in V \mid (\lambda,\alpha^\vee)\in\Z
\textrm{ for any } \alpha\in\Phi\}$.
It is generated by the {\it fundamental weights\/}
$\omega_1,\dots,\omega_r$ that form the dual basis to the 
basis of simple coroots, i.e., $(\omega_i,\alpha_j^\vee)=\delta_{ij}$.
The set $\Lambda^+$ of {\it dominant weights\/} is given by
$\Lambda^+=\{\lambda\in\Lambda \mid (\lambda,\alpha^\vee)\geq 0
\textrm{ for any } \alpha\in\Phi^+\}$.
A dominant weight $\lambda$ is called {\it regular\/} if 
$(\lambda,\alpha^\vee)>0$ for any $\alpha\in\Phi^+$.
Let $\rho=\omega_1+\cdots+\omega_r=\frac{1}{2}\sum_{\alpha\in\Phi^+}\alpha$
be the minimal regular dominant weight.


\section{Schubert calculus}
\label{sec:schubert-calc}

In this section, we recall some classical results of 
Borel~\cite{Bor}, Chevalley~\cite{Chev}, Demazure~\cite{Dem},
and Bernstein-Gelfand-Gelfand~\cite{BGG}.
\smallskip

The {\it generalized flag variety\/} $G/B$ is a smooth complex projective
variety.  Let $H^*(G/B)=H^*(G/B,\Q)$ be the {\it cohomology ring\/} 
of $G/B$ with rational coefficients.  
Let $\Q[V^*]=\Sym(V)$ be the algebra of polynomials on 
the space $V^*$ with rational coefficients.
The action of the Weyl group $W$ on the space $V$ induces a $W$-action 
on the polynomial ring $\Q[V^*]$.   
According to {\it Borel's theorem\/}~\cite{Bor}, the cohomology of $G/B$
is canonically isomorphic%
\footnote{The isomorphism is given by
$c_1(\L_\lambda)\mapsto \lambda\pmod{\I_W}$,
where $c_1(\L_\lambda)$ is the first Chern class of
the line bundle $\L_\lambda= G\times_B \C_{-\lambda}$ over $G/B$,
for $\lambda\in\Lambda^+$.}
to the quotient of the polynomial ring:
\begin{equation}
H^*(G/B)\simeq \Q[V^*]/\I_W,
\label{eq:Borel}
\end{equation}
where $\I_W= \left<f\in \Q[V^*]^W\mid f(0)=0\right>$ 
is the ideal generated by $W$-invariant 
polynomials without constant term.
Let us identify the cohomology ring $H^*(G/B)$ with this quotient 
ring. 
For a polynomial $f\in \Q[V^*]$,
let $\bar f=f\pmod\I_W$ be its coset modulo $\I_W$, which we view as 
a class in the cohomology ring $H^*(G/B)$.

One can construct a linear basis of $H^*(G/B)$ using 
the following {\it divided difference operators\/}
(also known as the Bernstein-Gelfand-Gelfand operators).
For a root $\alpha\in\Phi$, let $A_\alpha:\Q[V^*] \to \Q[V^*]$
be the operator given by 
\begin{equation}
A_\alpha:f\mapsto \frac {f - s_\alpha(f)} {\alpha}.
\label{eq:BGG-operator}
\end{equation}
Notice that the polynomial $f-s_\alpha(f)$ is always divisible by $\alpha$.
The operators $A_\alpha$ commute with operators of multiplication by 
$W$-invariant polynomials.
Thus the $A_\alpha$ preserve the ideal $\I_W$ and induce operators acting
on $H^*(G/B)$, which we will denote by the same symbols $A_\alpha$.

Let $A_i=A_{\alpha_i}$, for $i=1,\dots,r$.
The operators $A_i$ satisfy the nilCoxeter relations 
$$
(A_i A_j)^{m_{ij}}=1\qquad\textrm{and}\qquad (A_i)^2=0.
$$ 
For a reduced decomposition $w=s_{i_1}\cdots s_{i_l}\in W$,
define $A_w=A_{i_1}\cdots A_{i_l}$.  The operator $A_w$ depends only on 
$w\in W$ and does not depend on a choice of reduced decomposition.

Let us define the {\it Schubert classes\/} $\sigma_w\in H^*(G/B)$, $w\in W$, 
by
$$
\begin{array}{l}
\displaystyle
\sigma_{w_\circ}=|W|^{-1}\, \prod_{\alpha\in\Phi^+} \alpha 
\pmod {\I_W},\quad \textrm{for the longest element }w_\circ\in W;\\[.15in]
\displaystyle
\sigma_w=A_{w^{-1}w_\circ}(\sigma_{w_\circ}),
\quad\textrm{for any }w\in W.
\end{array}
$$

The classes $\sigma_w$ have the following geometrical meaning.
Let $X_w = \overline{BwB/B}$, $w\in W$, be the {\it Schubert varieties\/} 
in $G/B$.  According to Bernstein-Gelfand-Gelfand~\cite{BGG} and 
Demazure~\cite{Dem}, 
$\sigma_w = [X_{w_\circ w}]\in H^{2\ell(w)}(G/B)$ 
are the cohomology classes of the Schubert varieties.
They form a linear basis of the cohomology ring $H^*(G/B)$.
In the basis of Schubert classes, the divided difference operators 
can be expressed, as follows (see~\cite{BGG}):
\begin{equation}
A_i(\sigma_{w}) = 
\left\{
\begin{array}{cl}
\sigma_{ws_i} &\textrm{if } \ell(ws_i)=\ell(w)-1,\\
0 &\textrm{if } \ell(ws_i)=\ell(w)+1.
\end{array}
\right.
\label{eq:BGG-formula}
\end{equation}

\begin{remark}
There are many possible choices for polynomial
representatives of the Schubert classes.
In type $A_{n-1}$, Lascoux and Sch\"utzenberger~\cite{LS} introduced
the polynomial representatives, called the {\it Schubert polynomials,} 
obtained from the monomial
$x_1^{n-1}x_2^{n-2}\cdots x_{n-1}$ by applying the divided difference
operators.
Here $x_1,\dots,x_n$ are the coordinates 
in the standard presentation for type $A_{n-1}$ roots
$\alpha_{ij} = x_i - x_j$  (see~\cite{Hum}).
Schubert polynomials 
have many nice combinatorial properties;
see Section~\ref{sec:A-Schubert-poly} below.
\end{remark}


For $\sigma\in H^*(G/B)$, let $\<\sigma\>=\int_{G/B}\sigma$ 
be the coefficient of the top class $\sigma_{w_\circ}$ 
in the expansion of $\sigma$ in the Schubert classes.
Then $\<\sigma\cdot \theta\>$ is the 
{\it Poincar\'e pairing\/} on $H^*(G/B)$.  In the basis of Schubert 
classes the Poincar\'e pairing is given by
\begin{equation}
\<\sigma_u\cdot \sigma_{w}\> = \delta_{u,\,w_\circ w}.
\label{eq:poincare}
\end{equation}

The {\it generalized Littlewood-Richardson coefficients\/} $c_{u,v}^w$,
are given by
$$
\sigma_u\cdot\sigma_v = \sum_{w\in W} c_{u,v}^w\, \sigma_{w}, \quad
\textrm{for }u,v\in W. 
$$
Let $c_{u,v,w} = \<\sigma_u\cdot \sigma_v\cdot \sigma_{w}\>$ be the 
triple intersection number of Schubert varieties.  Then, according 
to~(\ref{eq:poincare}), we have $c_{u,v}^w=c_{u,v,w_\circ w}$.

\begin{figure}[ht]
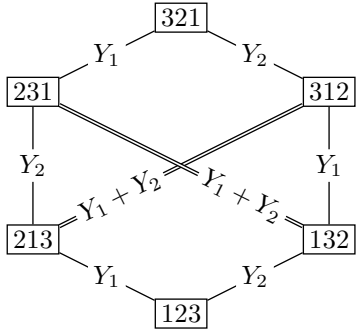

\pspicture(0,0)(0,90)
\rput(0,0){\rnode{A}{\psframebox{123}}}
\rput(-40,20){\rnode{B}{\psframebox{213}}}
\rput(40,20){\rnode{C}{\psframebox{132}}}
\rput(-40,60){\rnode{D}{\psframebox{231}}}
\rput(40,60){\rnode{E}{\psframebox{312}}}
\rput(0,80){\rnode{F}{\psframebox{321}}}
\ncline{-}{A}{B}
\mput*{$Y_1$}
\ncline{-}{A}{C}
\mput*{$Y_2$}
\ncline{-}{B}{D}
\mput*{$Y_2$}
\ncline[doubleline=true]{-}{B}{E}
\lput*{:U}(.25){$Y_1+Y_2$}
\ncline[doubleline=true]{-}{D}{C}
\lput*{:U}(.75){$Y_1+Y_2$}
\ncline{-}{C}{E}
\mput*{$Y_1$}
\ncline{-}{D}{F}
\mput*{$Y_1$}
\ncline{-}{E}{F}
\mput*{$Y_2$}
\endpspicture
\caption{The Bruhat order on $S_3$ marked with 
the Chevalley multiplicities.}
\label{fig:1}
\end{figure}

For a linear form $y\in V\subset \Q[V^*]$, let 
$\bar y \in H^*(G/B)$ be its coset%
\footnote{Equivalently, $\bar y = c_1(\L_\lambda)$, 
if $y=\lambda$ is in the weight lattice $\Lambda$.}
modulo $\I_W$.
{\it Chevalley's formula}~\cite{Chev} gives the following rule for the 
product of a Schubert class $\sigma_w$, $w\in W$, with $\bar y$:
\begin{equation}
\bar y\cdot \sigma_w=\sum
(y,\alpha^\vee)\,\sigma_{w s_\alpha},
\label{eq:Chevalley}
\end{equation}
where the sum is over all roots $\alpha\in\Phi^+$ such that
$\ell(w\,s_\alpha)=\ell(w)+1$, i.e., the sum is over all elements in $W$
that cover $w$ in the Bruhat order.
The coefficients $(y,\alpha^\vee)$, which are associated
to edges in the Hasse diagram of the Bruhat order, are called
the {\it Chevalley multiplicities}.
Figure~\ref{fig:1} shows the Bruhat order on the symmetric
group $W=S_3$ with edges of the Hasse diagram marked 
by the Chevalley multiplicities,
where $Y_1=(y,\alpha_1^\vee)$ and $Y_2=(y,\alpha_2^\vee)$.

We have, $\sigma_{\id}=[G/B]=1$.  Chevalley's formula implies that 
$\sigma_{s_i}=\bar \omega_i$ (the coset of the fundamental weight
$\omega_i$).

\psset{unit=1.4pt} 
\psset{linewidth=.5pt}
\psset{framesep=2pt}

\section{Degrees of Schubert varieties}
\label{sec:deg-schub-var}

For $y\in V$, let $m(u\lessdot u s_\alpha) = (y,\alpha^\vee)$ 
denote the Chevalley multiplicity of a covering relation 
$u \lessdot u \,s_\alpha$ in the Bruhat order
on the Weyl group $W$. 
Let us define the {\it weight\/} $m_C=m_C(y)$ of a saturated chain 
$C=(u_0\lessdot u_1 \lessdot u_2\lessdot \cdots \lessdot u_l)$
in  the Bruhat order as the product of Chevalley multiplicities:
$$
m_C(y) = \prod_{i=1}^l m(u_{i-1}\lessdot u_i).
$$ 
Then the weight $m_C\in\Q[V]$ is a polynomial function 
of $y\in V$.

For two Weyl group elements $u,w\in W$, $u\leq w$,  
let us define the polynomial $\D_{u,w}(y)\in\Q[V]$ as the sum
\begin{equation}
\D_{u,w}(y)= 
\frac{1}{(\ell(w)-\ell(u))!}\, \sum_C m_C(y)
\label{eq:D-uw-def}
\end{equation}
over all saturated chains
$C=(u_0\lessdot u_1 \lessdot u_2\lessdot \cdots \lessdot u_l)$
in the Bruhat order from $u_0=u$ to $u_l=w$. 
In particular, $\D_{w,w}=1$.
Let $\D_w=\D_{\id,w}$.
It is clear from the definition that 
$\D_w$ is a homogeneous polynomial of degree $\ell(w)$ and
$\D_{u,w}$ is homogeneous of degree $\ell(w)-\ell(u)$. 

\begin{example}
For $W=S_3$, we have
$\D_{\id, 231} = \frac{1}{2}(Y_1Y_2 + Y_2(Y_1+Y_2))$ and
$\D_{132, 321} = \frac{1}{2}((Y_1+Y_2)Y_1 + Y_1 Y_2)$,
where $Y_1=(y,\alpha_1^\vee)$ and $Y_2=(y,\alpha_2^\vee)$
(see Figure~\ref{fig:1}).
\end{example}

According to Chevalley's formula~(\ref{eq:Chevalley}), 
the values of the polynomials $\D_{u,w}(y)$ 
are the expansion coefficients in the following product in 
the cohomology ring $H^*(G/B)$:
\begin{equation}
[e^y]\cdot \sigma_u = \sum_{w\in W} 
\D_{u,w}(y)\cdot \sigma_w, 
\textrm{ for any } y\in V,
\label{eq:e-lambda-D-uw}
\end{equation}
where $[e^y]:=1+\bar y + \bar y^2/2!+
\bar y^3/3!+\cdots\in H^*(G/B)$.
Note that $[e^y]$ involves only finitely many nonzero summands, 
because $H^{k}(G/B)=0$, for sufficiently large $k$.
Equation~(\ref{eq:e-lambda-D-uw}) is actually equivalent to
definition~(\ref{eq:D-uw-def}) of the polynomials $\D_{u,w}$.

The values of the polynomials $\D_w(\lambda)$ at dominant weights 
$\lambda\in\Lambda^+$
have the following natural geometric interpretation. 
For $\lambda\in\Lambda^+$,
let $V_\lambda$ be the {\it irreducible representation\/} of the Lie group
$G$ with the highest weight $\lambda$, and let
$v_\lambda\in V_\lambda$ be a highest weight vector. 
Let $e:G/B\to\mathbb{P}(V_\lambda)$ be the map given by 
$gB\mapsto g(v_\lambda)$, for $g\in G$.  If the weight $\lambda$ is regular,
then $e$ is a projective embedding $G/B\hookrightarrow\mathbb{P}(V_\lambda)$.
Let $w\in W$ be an element of length $l=\ell(w)$.
Let us define the {\it $\lambda$-degree\/} $\deg_\lambda(X_w)$ 
of the Schubert variety $X_w\subset G/B$ 
as the number of points in the intersection of $e(X_w)$ with a generic
linear subspace in $\mathbb{P}(V_\lambda)$ of complex codimension $l$.
The pull-back of the class of a hyperplane in
$H^2(\mathbb{P}(V_\lambda))$ is 
$\bar\lambda=c_1(\L_\lambda)\in H^2(G/B)$.  
Then the $\lambda$-degree of $X_w$
is equal to the Poincar\'e pairing 
$\deg_\lambda(X_w) = \left<[X_w]\cdot \bar \lambda^l\right>$. 
In other words, $\deg_\lambda(X_w)$ equals the coefficient of 
the Schubert class $\sigma_w$, which is Poincar\'e dual to
$[X_w]=\sigma_{w_\circ w}$, in the expansion of $\bar \lambda^l$ 
in the basis of Schubert classes.
Chevalley's formula~(\ref{eq:Chevalley})
implies the following well-known statement; see, e.g.,~\cite{BL}.

\begin{proposition}  For $w\in W$ and $\lambda\in\Lambda^+$,
the $\lambda$-degree $\deg_\lambda(X_w)$ of the Schubert variety $X_w$ 
is equal to the sum $\sum m_C(\lambda)$
over saturated chains $C$ in the Bruhat order from $\id$ to $w$.
Equivalently,
$$
\deg_\lambda(X_w)=\ell(w)!\cdot \D_w(\lambda).
$$
\end{proposition}

If $\lambda=\rho$, we will call $\deg(X_w)=\deg_\rho(X_w)$
simply the {\it degree\/} of $X_w$.

\section{Harmonic polynomials}
\label{sec:D-pairing}

We discuss harmonic polynomials and the natural pairing on polynomials 
defined in terms of partial derivatives.   Constructions in this section 
are essentially well-known; cf.~Bergeron-Garsia~\cite{BG}.
\smallskip


The space of polynomials $\Q[V]$ is the graded dual to $\Q[V^*]$, i.e.,
the corresponding finite-dimensional graded components are dual
to each other.


Let us pick a basis $v_1,\dots,v_r$ in $V$, and let $v_1^*,\dots,v_r^*$ 
be the dual basis in $V^*$.  For $f\in\Q[V^*]$ and $g\in\Q[V]$,
let $f(x_1,\dots,x_r)=f(x_1v_1^*+\cdots+x_r v_r^*)$
and $g(y_1,\dots,y_r)=g(y_1 v_1+\cdots+y_r v_r)$
be polynomials in the variables $x_1,\dots,x_r$ and 
$y_1,\dots,y_r$, correspondingly.
For each $f\in\Q[V^*]$, let us define the differential operator $f(\ddy)$ 
that acts on the polynomial ring $\Q[V]$ by
$$
f(\ddy): g(y_1,\dots,y_r) \longmapsto
f(\ddy_1,\dots,\ddy_r)\cdot g(y_1,\dots,y_r),
$$
where $\ddy_i$ denotes the partial derivative with respect to $y_i$.
The operator $f(\ddy)$ can also be described without coordinates as follows.
Let $d_v:\Q[V]\to\Q[V]$ be the differentiation operator 
in the direction of a vector $v\in V$ given by
\begin{equation}
d_v: g(y)\mapsto \left.\frac{d}{d\,t}
\,g(y+t\,v)
\right|_{t=0}.
\label{eq:diff-mu}
\end{equation}
The linear map $v\mapsto d_v$ extends to the homomorphism $f\mapsto d_f$
from the polynomial ring $\Q[V]=\Sym(V^*)$ to the ring of operators 
on $\Q[V^*]$.  Then $d_f=f(\ddy)$ 

One can extend the usual pairing between $V$ and $V^*$
to the following pairing between the spaces $\Q[V^*]$ and $\Q[V]$. 
For $f\in\Q[V^*]$ and $g\in\Q[V]$, let us 
define the {\it D-pairing\/} $\(f,g\)_D$ by
$$
\(f,g\)_D =\CT\(f(\ddy)\cdot g(y)\)=
\CT\(g(\ddx)\cdot f(x)\),
$$
where the notation
$\CT$ means taking the constant term of a polynomial.

A graded basis of a polynomial ring is a basis that consists 
of homogeneous polynomials.
Let us say that a graded $\Q$-basis $\{f_u\}_{u\in U}$ in $\Q[V^*]$
is {\it D-dual\/} to 
a graded $\Q$-basis $\{g_u\}_{u\in U}$ in $\Q[V]$
if $(f_u,g_v)_D=\delta_{u,v}$, for any $u,v\in U$.

\begin{example}
Let $x^a = x_1^{a_1} \cdots x_r^{a_r}$ and 
$y^{(a)} = \frac{y_1^{a_1}}{a_1!} \cdots 
\frac{y_r^{a_r}}{a_r!}$, for $a=(a_1,\dots,a_r)$.
Then the monomial basis $\{x^a\}$ of $\Q[V^*]$ is D-dual to the basis
$\{y^{(a)}\}$ of $\Q[V]$.
\end{example}

This example shows that the D-pairing gives a non-degenerate pairing
of corresponding graded components of $\Q[V^*]$ and $\Q[V]$ and 
vanishes on different graded components.
Thus, for a graded basis in $\Q[V^*]$, there exists a unique 
D-dual graded basis in $\Q[V]$ and vice versa.

For a graded space $A=A^0\oplus A^1\oplus A^2\oplus \cdots$,
let $A_\infty$ be the space of formal series $a_0+a_1+a_2+\cdots$,
where $a_i\in A^i$.
For example, $\Q[V]_\infty = \Q[[V]]$ is the ring of formal power series.
The exponential $e^{(x,y)}= e^{x_1 y_1+\cdots + x_r y_r}$
given by its Taylor series can 
be regarded as an element of $\Q[[V^*,V]]:=\Q[[V^*]]\otimes \Q[[V]]$,
where $(x,y)$ is the standard pairing between $x\in V^*$ 
and $y\in V$.

\begin{proposition}
Let $\{f_u\}_{u\in U}$ be a graded basis for $\Q[V^*]$, 
and let $\{g_u\}_{u\in U}$ be a collection of formal power series in 
$\Q[[V]]$ labeled by the same set $U$.
Then the following two conditions are equivalent:
\begin{enumerate}
\item  
The $g_u$ are the homogeneous polynomials in $\Q[V]$ that 
form the D-dual basis to $\{f_u\}$.
\item
The equality 
$e^{(x,y)} = \sum_{u\in U} f_u(x)\cdot g_u(y)$
holds identically in the ring of formal power series $\Q[[V^*,V]]$.
\end{enumerate}
\label{prop:D-dual}
\end{proposition}

\begin{proof}
For $f\in\Q[V^*]$, the action of the differential operator $f(\ddy)$ on 
polynomials extends to the action on the ring of formal power series 
$\Q[[V]]$ and on $\Q[[V^*,V]]$.
The D-pairing $(f,g)_D$ makes sense for any $f\in\Q[V^*]$ and $g\in\Q[[V]]$.
Let $C=\sum_{u\in U} f_u(x)\cdot g_u(y)\in \Q[[V^*,V]]$.   
Then $\CT\(f_u(\ddy)\cdot C\) = \sum_{v\in U} 
(f_u,g_v)_D\, f_v(x)$, for any $u\in U$.

Condition (1) is equivalent to the condition $\CT\(f(\ddy)\cdot C\) = f(x)$, 
for any basis element $f=f_u$ of $\Q[V^*]$.  The latter condition is 
equivalent to condition (2), which says that $C=e^{(x,y)}$.
Indeed, the only element $E\in \Q[[V^*,V]]$ that satisfies
$\CT\(f(\ddy)\cdot E\)=f(x)$, for any
$f\in \Q[V^*]$,  is the exponent $E=e^{(x,y)}$.
%
\end{proof}

Let $I\subseteq \Q[V^*]$ be a graded ideal. 
Define the space of {\it $I$-harmonic polynomials\/} as
$$
\H_I = 
\{g\in\Q[V]\mid f(\ddy)\cdot g(y) = 0, \textrm{ for any } f\in I\}.
$$

\begin{lemma}  
The space $\H_I\subseteq \Q[V]$ 
is the orthogonal subspace to $I\subseteq\Q[V^*]$ with respect to the D-pairing.
Thus $\H_I$ is the graded dual to the quotient space $\Q[V^*]/I$.
\end{lemma}

\begin{proof}
The ideal $I$ is orthogonal to 
$I^\perp:=\{g\mid (f,g)_D=0, \textrm{ for any } f\in I\}$.
Clearly, $\H_I\subseteq I^\perp$.  On the other hand, if 
$(f,g)_D=\CT\(f(\ddy)\cdot g(y)\)=0$, for any $f\in I$,
then $f(\ddy)\cdot g(y)=0$, for any $f\in I$, because $I$
in an ideal.  Thus $\H_I=I^\perp$.
\end{proof}

Let $\bar f:=f\pmod I$ denote the coset of a polynomial $f\in\Q[V^*]$ modulo the ideal
$I$.  For $g\in \H_I$, the differentiation 
$\bar f(\ddy)\cdot g:=f(\ddy)\cdot g$ does not depend 
on the choice of a polynomial representative $f$ of the coset $\bar f$.
Thus we have correctly defined a D-pairing $(\bar f, g)_D:=(f, g)_D$ between 
the spaces $\Q[V^*]/I$ and $\H_I$.  Let us say that a graded basis 
$\{\bar f_u\}_{u\in U}$ of $\Q[V^*]/I$ and 
a graded basis $\{g_u\}_{u\in U}$ of $\H_I$ are {\it D-dual\/}
if $(\bar f_u, g_v)_D=\delta_{u,v}$, for any $u,v\in U$.

\begin{proposition}
Let $\{\bar f_u\}_{u\in U}$ be a graded basis of $\Q[V^*]/I$, 
and let $\{g_u\}_{u\in U}$ be a collection of 
of formal power series in $\Q[[V]]$ labeled by the same set $U$.
Then the following two conditions are equivalent:
\begin{enumerate}
\item  
The $g_u$ are the polynomials that form the graded
basis of $\H_I$ such that
the bases $\{\bar f_u\}_{u\in U}$ and $\{g_u\}_{u\in U}$
are D-dual.
\item
The equality 
$e^{(x,y)} = \sum_{u\in U} f_u(x)\cdot g_u(y)$
  modulo $I_\infty\otimes \Q[[V]]$
holds identically.
\end{enumerate}
\label{prop:exlambdamod}
\end{proposition}

\begin{proof}
Let us augment the set $\{f_u\}_{u\in U}$ by
a graded $\Q$-basis $\{f_{u}\}_{u\in U'}$ of the ideal $I$.
Then $\{f_u\}_{u\in U\cup U'}$  is a graded basis of $\Q[V^*]$.
A collection $\{g_u\}_{u\in U}$ is the basis of $\H_I$ that is D-dual to 
$\{\bar f_u\}_{u\in U}$ if and only if there are elements 
$g_{u}\in\Q[V]$, for $u\in U'$,
such that $\{f_u\}_{u\in U\cup U'}$ and 
$\{g_u\}_{u\in U\cup U'}$ are D-dual bases
of $\Q[V^*]$ and $\Q[V]$, correspondingly.
The claim now follows from Proposition~\ref{prop:D-dual}.
\end{proof}

The {\it product map\/} $M:\Q[V^*]/I\otimes\Q[V^*]/I\to\Q[V^*]/I$ is given by
$M:\bar f\otimes \bar g\mapsto \bar f\cdot \bar g$.
Let us define {\it coproduct map\/} $\Delta:\H_I\to\H_I\otimes \H_I$
as the D-dual map to $M$.
For $h\in\Q[V]$, the polynomial $h(y+z)$ of the sum of two vector
variables $y,z\in V$ can be regarded as an element of $\Q[V]\otimes\Q[V]$.

\begin{proposition}  The coproduct map $\Delta:\H_I\to\H_I\otimes \H_I$ 
is given by
$$
\Delta:g(y)\mapsto g(y+z),
$$
for any $g\in\H_I$.
\label{prop:Delta-g}
\end{proposition}

\begin{proof} Let $\{\bar f_u\}_{u\in U}$ be a graded basis in $\Q[V^*]/I$
and let $\{g_u\}_{u\in U}$ be its D-dual basis in $\H_I$.
We need to show that the two expressions
$$
\bar f_u(x)\cdot \bar f_v(x) = \sum_{w\in U} a_{u,v}^w\,\bar f_w(x)
\quad\textrm{and}\quad
\bar g_w(y+z)= \sum_{u,v\in U} b_{u,v}^w\,g_u(y)\cdot g_v(z)
$$
have the same coefficients $a_{u,v}^w= b_{u,v}^w$.
Here $x\in V^*$ and $y,z\in V$.
Indeed, according to  
Proposition~\ref{prop:exlambdamod}, we have
$$
\begin{array}{l}
\displaystyle
\sum_{u,v,w} a_{u,v}^w \,\bar f_w(x)\cdot g_u(y)\cdot g_v(z)=
\left(\sum_{u} \bar f_u(x)\cdot g_u(y)\right)\cdot
\left(\sum_{v} \bar f_v(x)\cdot g_v(z)\right)=\\[.2in]
\displaystyle
=e^{(x,y)}\, e^{(x,z)} = e^{(x,y+z)} =
\sum_{w} \bar f_w(x)\cdot g_w(y+z)
=\sum_{u,v,w} b_{u,v}^w\,\bar f_w(x)
\cdot g_u(y)\cdot g_v(z)
\end{array}
$$
in the space $(\Q[V^*]/I\otimes \Q[V]\otimes \Q[V])_\infty$.
This implies that $a_{u,v}^w=b_{u,v}^w$, for any $u,v,w\in U$.
\end{proof}

In what follows, we will assume 
and $I=\I_W\subset \Q[V^*]$ is the ideal generated by 
$W$-invariant polynomials without
constant term,  
and $\Q[V^*]/I= H^*(G/B)$ is the cohomology ring of $G/B$.
Let $\H_W=\H_{\I_W}\subset\Q[V]$ be its dual space with respect 
to the D-pairing.  We will call $\H_W$ the {\it space 
of $W$-harmonic polynomials\/}
and call its elements {\it $W$-harmonic polynomials\/} in $\Q[V]$.

\section{Expressions for polynomials $\D_{u,w}$}
\label{sec:expr-for-D}

In this section, we give two different expressions for 
the polynomials $\D_{u,w}$ and derive several corollaries.
\smallskip

Formula~(\ref{eq:e-lambda-D-uw}), for $u=\id$, 
and Proposition~\ref{prop:exlambdamod} imply the following statement.

\begin{corollary}
{\rm (cf.~Bernstein-Gelfand-Gelfand~\cite[Theorem~3.13]{BGG})}
The collection of polynomials $\D_w$, $w\in W$, forms a linear basis 
of the space $\H_W\subset \Q[V]$ of $W$-harmonic polynomials.  
This basis is D-dual to the basis $\{\sigma_w\}_{w\in W}$ of 
Schubert classes in $H^*(G/B)$.
\label{cor:D-dual-sigma}
\end{corollary}

This basis of $W$-harmonic polynomials appeared in 
Bernstein-Gelfand-Gelfand \cite[Theorem~3.13]{BGG} 
(in somewhat disguised form) and more recently in Kriloff-Ram 
\cite[Sect.~2.2]{KR}; see Remark~\ref{rem:KR} below.

By the definition, the polynomial $\D_{u,w}$ is given by a sum over 
saturated chains in the Bruhat order.  However, this expression involves
many summands and is difficult to handle.
The following theorem given a more explicit formula for $\D_{u,w}$.

Let $\sigma_w(\ddy)$ be the differential operator 
on the space of $W$-harmonic polynomials $\H_W$ given by
$\sigma_w(\ddy): g(y) \mapsto \S_w(\ddy)\cdot g(y)$, 
where $\S_w\in\Q[V^*]$ is any polynomial representative of 
the Schubert class $\sigma_w$.
According to Section~\ref{sec:D-pairing}, $\sigma_w(\ddy)$ does not
depend on the choice of a polynomial representative $\S_w$.

\begin{theorem}
For any $w\in W$, we have
$$
\D_{u,w}(y) = \sigma_u(\ddy)\,\sigma_{w_\circ w}(\ddy)\cdot \D_{w_\circ}(y).
$$
In particular, all polynomials $\D_{u,w}$ are $W$-harmonic.
\label{th:D=S}
\end{theorem}  


\begin{proof}
According to~(\ref{eq:e-lambda-D-uw}), 
we have $\D_{u,w}(\lambda)=\<[e^\lambda]\cdot 
\sigma_u\cdot \sigma_{w_\circ w}\>$, for any weight $\lambda\in\Lambda$.
Thus the $W$-harmonic polynomial $\D_{u,w}$ is uniquely determined by the identity
$(\sigma,\D_{u,w})_D = \<\sigma\cdot \sigma_u\cdot \sigma_{w_\circ w}\>$,
for any $\sigma\in H^*(G/B)$.
Let us show that the same identity holds for the $W$-harmonic polynomial
$\tilde\D_{u,w}(y)=\sigma_u(\ddy)\,\sigma_{w_\circ w}(\ddy)\cdot 
\D_{w_\circ}(y)$.
Indeed, $(\sigma,\tilde\D_{u,w})_D$ equals
$$
\CT\(\sigma(\ddy)\cdot \sigma_u(\ddy)\cdot \sigma_{w_\circ w}(\ddy)\cdot 
\D_{w_\circ}(y)\)
=(\sigma\cdot \sigma_u \cdot \sigma_{w_\circ w},\D_{w_\circ})_D.
$$
Since $\{\D_w\}_{w\in W}$ is the D-dual basis to $\{\sigma_w\}_{w\in W}$,
the last expression is equal to triple intersection number
$\<\sigma\cdot \sigma_u\cdot \sigma_{w_\circ w}\>$, as needed.
\end{proof}

Corollary~\ref{cor:D-w0} below gives a simple multiplicative 
Vandermonde-like expression
for $\D_{w_\circ}$.  Theorem~\ref{th:D=S}, together with 
this expression, gives an explicit ``top-to-bottom'' 
differential formula for the $W$-harmonic polynomials $\D_w$.
Let us give an alternative ``bottom-to-top'' integral formula for 
these polynomials.

For $\alpha\in\Phi$, let $I_\alpha$ be the operator that acts on polynomials
$g\in\Q[V]$ by 
\begin{equation}
I_i:g(y)\mapsto \int_{0}^{(y,\alpha^\vee)}
g(y-\alpha t)\,dt.
\label{eq:Ii}
\end{equation}
In other words, the operator $I_\alpha$ integrates a polynomial $g$ 
on the line interval $[y,s_\alpha(y)]\subset V$.  Clearly, this operator
increases the degree of polynomials by 1.

Recall that $A_\alpha:\Q[V^*]\to\Q[V^*]$ is the BGG operator given 
by~(\ref{eq:BGG-operator}).

\begin{lemma}  For $\alpha\in\Phi$, the operator $I_\alpha$ is adjoint
to the operator $A_\alpha$ with respect to the D-pairing.
In other words,
\begin{equation}
(f,I_\alpha(g))_D = (A_\alpha(f),g)_D,
\label{eq:I-conjugate-A}
\end{equation}
for any polynomials $f\in\Q[V^*]$ and $g\in\Q[V]$.
\label{lem:I-conjugate-A}
\end{lemma}

\begin{proof} Let us pick a basis $v_1,\dots,v_r$ in $V$
and its dual basis $v_1^*,\dots,v_r^*$ in $V^*$ 
such that $v_1=\alpha$ and $(v_i,\alpha)=0$, for $i=2,\dots,r$.
Let $f(x_1,\dots,x_r)=f(x_1v_1^*+\cdots+x_r v_r^*)$ and
$g(y_1,\dots,y_r)= 
g(y_1 v_1+\cdots + y_r v_r)$, for $f\in\Q[V^*]$
and $g\in\Q[V]$.  In these coordinates, the operators 
$A_\alpha$ and $I_\alpha$ can be written as
$$
\begin{array}{l}
\displaystyle
A_\alpha:f(x_1,\dots,x_r)\mapsto 
\frac{f(x_1,x_2,\dots,x_r)-f(-x_1,x_2,\dots,x_r)}{x_1}\\[.2in]
\displaystyle
I_\alpha:g(y_1,\dots,y_r)\mapsto
\int_{-y_1}^{y_1} g(t,y_2,\dots,y_r)\,dt.
\end{array}
$$
These operators are linear over $\Q[x_2,\dots,x_r]$ and 
$\Q[y_2,\dots,y_r]$, correspondingly.  
It is enough to verify identity~(\ref{eq:I-conjugate-A}) for 
$f=x_1^{m+1}$ and $g=y_1^{m}$.
For these polynomials, we have $A_\alpha(f) = 2x_1^m$,
$I_\alpha(g) = \frac 2{m+1} y_1^{m+1}$, if $m$ is even;
and $A_\alpha(f) = 0$, $I_\alpha(g)=0$, if $m$ is odd.
Thus $(f,I_\alpha(g))_D = (A_\alpha(f),g)_D$
in both cases.
\end{proof}

Let $I_i=I_{\alpha_i}$, for $i=1,\dots,r$.

\begin{corollary}  The operators $I_i$ satisfy the 
nilCoxeter relations $(I_i I_j)^{m_{ij}}=1$ and $(I_i)^2=0$.
Also, if $I_\alpha(g) = 0$, then $g$ is an anti-symmetric polynomial 
with respect to the reflection $s_\alpha$, and thus, 
$g$ is divisible by the linear form $(y,\alpha^\vee)\in\Q[V]$.
\label{cor:divisible-by-linear-form}
\end{corollary}

\begin{proof}  The first claim follows from the fact that 
the BGG operators $A_i$ satisfy the nilCoxeter relations.
The second claim is clear from the formula for $I_\alpha$
given in the proof of Lemma~\ref{lem:I-conjugate-A}.
\end{proof}

For a reduced decomposition $w=s_{i_1}\cdots s_{i_l}$, let us 
define $I_w= I_{i_1}\cdots I_{s_l}$.  The operator $I_w$ depends only 
on $w$ and does not depend on the choice of reduced decomposition.
Lemma~\ref{lem:I-conjugate-A} implies that the operator 
$A_w:\Q[V^*]\to\Q[V^*]$ is
adjoint to the operator $I_{w^{-1}}:\Q[V]\to\Q[V]$ with respect
to the D-pairing.

\begin{theorem}
{\rm (cf.\ Bernstein-Gelfand-Gelfand~\cite[Theorem~3.12]{BGG})}
For any $w\in W$ and $i=1,\dots,r$, we have
$$
I_i\cdot \D_w =
\left\{
\begin{array}{cl}
\D_{ws_i} &\textrm{if } \ell(ws_i)=\ell(w)+1,\\
0 &\textrm{if } \ell(ws_i)=\ell(w)-1.
\end{array}
\right.
$$
Thus the polynomials $\D_w$ are given by
$$
\D_w = I_{w^{-1}}(1).
$$
\label{th:integration-formula}
\end{theorem}

\begin{proof}
Follows from Bernstein-Gelfand-Gelfand formula~(\ref{eq:BGG-formula}),
Corollary~\ref{cor:D-dual-sigma}, and Lemma~\ref{lem:I-conjugate-A}.
\end{proof}

\begin{remark}  
Theorem~\ref{th:integration-formula} is essentially contained 
in~\cite{BGG}.
However, Bernstein-Gelfand-Gelfand treated the $\D_w$ not as 
(harmonic) polynomials but as linear functionals on $\Q[V^*]/\I_W$
obtained from $\mathrm{Id}$ by appying the operators  
adjoint to the divided difference operators operators (with respect to 
the natural pairing between a linear space and its dual).
It is immediate that these functionals form a basis 
in  $(\Q[V^*]/\I_W)^*\simeq \H_W$;
see~\cite[Theorem~3.13]{BGG} and~\cite[Sect.~2.2]{KR}.
Note that there are several other constructions of bases of $\H_W$;
see, e.g., Hulsurkar~\cite{Hul}.
\label{rem:KR}
\end{remark}

In the next section we show that Duan's recent result~\cite{Du}
about degrees of Schubert varieties easily follows from
Theorem~\ref{th:integration-formula}.
Note that this integral expression for the polynomials $\D_w$ can 
be formulated in the general Kac-Moody setup.  Indeed, unlike 
the previous expression given by Theorem~\ref{th:D=S}, it does not use 
the longest Weyl group element $w_\circ$, 
which exists in finite types only.

For $I\subseteq\{1,\dots,r\}$, let $W_I$ be the parabolic 
subgroup in $W$ generated by $s_i$, $i\in I$.
Let $\Phi_I^+=\{\alpha\in\Phi^+\mid s_\alpha\in W_I\}$.

\begin{proposition}  Let $w\in W$. Let $I=
\{i\mid \ell(w s_i)<\ell(w)\}$ be the descent set of $w$.
Then the polynomial $\D_{w}(y)$ 
is divisible by the product 
$\prod_{\alpha\in \Phi_I^+} (y,\alpha^\vee)\in \Q[V]$.
\label{prop:divisible-by-product}
\end{proposition}

\begin{proof}  
According to Corollary~\ref{cor:divisible-by-linear-form},
it is enough to check that
$I_\alpha(\D_w)=0$, for any $\alpha\in\Phi_I^+$.
We have $I_\alpha(\D_w)=I_\alpha I_{w^{-1}}(1)$.
The operator $I_\alpha I_{w^{-1}}$ is adjoint to $A_w A_\alpha$
with respect to the D-pairing.
Let us show that $A_w A_\alpha=0$, identically.
Notice that $s_i A_\alpha = A_{s_i(\alpha)} s_i$,
where $s_i$ is regarded as an operator on the polynomial ring $\Q[V^*]$.
Also $A_i=s_iA_i=-A_i s_i$.
Thus, for any $i$ in the descent set $I$, we can write 
$$
A_w A_\alpha = A_{w'} A_i A_\alpha = - A_{w'} A_i s_i A_\alpha
= - A_{w'} A_i A_{s_i(\alpha)} s_i
= -A_w A_{s_i(\alpha)} s_i,
$$
where $w'=w s_i$.  Since $s_\alpha\in W_I$,
there is a sequence $i_1,\dots,i_l\in I$ and $j\in I$ such that
$s_{i_1} \cdots s_{i_l}(\alpha) = \alpha_j$. 
Thus 
$$
A_w A_\alpha = \pm A_w A_j s_{i_1}\cdots s_{i_l}=
\pm A_{w''} A_j A_j s_{i_1}\cdots s_{i_l}=0,
$$
as needed.
\end{proof}

\begin{corollary}  Fix $I\subseteq\{1,\dots,r\}$.
Let $w_I$ be the longest element in the parabolic subgroup $W_I$.
Then
$$
\D_{w_I}(y) = \mathrm{Const}\cdot\prod_{\alpha\in\Phi_I^+}
(y,\alpha^\vee),
$$
where $\mathrm{Const}\in\Q$.
\label{cor:D=const.prod}
\end{corollary}

\begin{proof}
Proposition~\ref{prop:divisible-by-product} says that
the polynomial $\D_{w_I}(y)$ is divisible by the product 
$\prod_{\alpha\in\Phi_I^+}(y ,\alpha^\vee)$.
Since the degree of this polynomial equals
$$
\deg \D_{w_I} = \ell(w_I) = |\Phi_I^+| = \deg 
\prod_{\alpha\in\Phi_I^+}(y ,\alpha^\vee),
$$
we deduce the claim.
\end{proof}

In Section~\ref{sec:asymtotic} below,
we will give an alternative derivation for this multiplicative
expression for $\D_{w_I}$; see Corollary~\ref{cor:D-w0}.
We will show that the constant $\mathrm{Const}$ in 
Corollary~\ref{cor:D=const.prod} is given by the condition
$\D_{w_I}(\rho)=1$.

We can express the generalized Littlewood-Richardson coefficients
$c_{u,v}^w$ using the polynomials $\D_{u,w}$ in two different ways.

\begin{corollary}
For any $u\leq w$ in $W$, we have
$$
\D_{u,w}=\sum_{v\in W} c_{u,v}^w\, \D_v.
$$
\label{Duwcw}
\end{corollary}

The polynomials $\D_{u,w}$ extend the polynomials $\D_{v}$ in 
the same way as the skew Schur polynomials
extend the usual Schur polynomials.
Compare Corollary~\ref{Duwcw} with a similar formula for 
the skew Schubert polynomials of Lenart and Sottile~\cite{LeS}.

\begin{proof}  
Let us expand the $W$-harmonic polynomial $\D_{u,w}$ in the basis 
$\{\D_v\mid v\in W\}$.  According to Theorem~\ref{th:D=S}, 
the coefficient of $\D_v$ in this expansion is equal to the 
coefficient of $\sigma_{w_\circ v}$ in the expansion 
of the product $\sigma_u\cdot \sigma_{w_\circ w}$ in the Schubert classes.
This coefficient equals $c_{u,w_\circ w}^{w_\circ v} = c_{u,v,w_\circ w} =
c_{u,v}^w$.
\end{proof}

Proposition~\ref{prop:Delta-g} implies the following statement.

\begin{corollary}  For $w\in W$, we have the 
equality\footnote{Here $y+z$ denotes the usual sum of two vectors.
This notation should not be confused with the $\lambda$-ring
notation for symmetric functions, where $y+z$ means the union of two
sets of variables.}

$$
\D_w(y+z)=\sum_{u,v\in W} c_{u,v}^w\, \D_u(y)\cdot \D_v(z)
$$
of polynomials in $y,z\in V$.
\label{cor:coporoduct}
\end{corollary}

Compare Corollary~\ref{cor:coporoduct} with the 
coproduct formula~\cite[Eq.~(7.66)]{EC2} for Schur polynomials.

\section{Examples and Duan's formula}
\label{sec:examp-duan}

Let us calculate several polynomials $\D_w$ using
Theorem~\ref{th:integration-formula}.
Let $Y_1,\dots,Y_r$ be the generators of $\Q[V]$
given by $Y_i=(y,\alpha_i^\vee)$, and let $a_{ij} = 
(\alpha_i^\vee,\alpha_j)$ be the Cartan integers,
for $1\leq i,j\leq r$.
For a simple reflection $w= s_i$, we obtain
$$
\D_{s_i} = I_i(1)=\int_0^{(y,\alpha_i^\vee)} 1\cdot dt = (y,\alpha_i^\vee) =
Y_i.
$$
For $w=s_i s_j$, we obtain
$$
\begin{array}{l}
\displaystyle
\D_{s_is_j} = I_j I_i(1) = I_j(Y_i)=I_j((y_,\alpha_i^\vee))=
\int_0^{(y,\alpha_j^\vee)} (y-t\,\alpha_j,\,\alpha_i^\vee)\,dt =\\[.2in]
\displaystyle
\qquad\qquad=(y,\alpha_i^\vee)\int_0^{(y,\alpha_j^\vee)} dt - 
(\alpha_j,\alpha_i^\vee) 
\int_0^{(y,\alpha_j^\vee)} t\,dt 
= Y_i Y_j - a_{ij} \,\frac{Y_j^2}{2}.
\end{array}
$$

We can further iterate this procedure.
The following lemma is obtained immediately from the definition 
of $I_j$'s, as shown above.

\begin{lemma}
For any $1\leq i_1,\dots,i_n,j\leq r$ and $c_1,\dots,c_n\in\Z_{\geq 0}$,
the operator $I_j$ maps the monomial
$Y_{i_1}^{c_1}\cdots Y_{i_n}^{c_n}$ to
$I_j(Y_{i_1}^{c_1}\cdots Y_{i_n}^{c_n}) = $
$$
\sum_{k_1+\cdots + k_n=k}
(-1)^k\,
\binom {c_1}{k_1} \cdots \binom {c_n}{k_n} \,
a_{i_1\,j}^{k_1}\cdots a_{i_n\,j}^{k_n} \,
Y_{i_1}^{c_1-k_1}\cdots Y_{i_n}^{c_n-k_n} \,
\frac{Y_j^{k+1}}{k+1}\,,
$$
where the sum is over $k_1,\dots,k_n$
such that $\sum k_i = k$, 
 $0\leq k_i\leq c_i$, for $i=1,\dots, n$.
\label{lem:I_i(jjj)}
\end{lemma}

For example, for $w=s_i s_j s_k$, we obtain
$$
\begin{array}{l}
\D_{s_i s_j s_k}  = I_k I_j I_i(1) = 
I_k(Y_iY_j-a_{ij}\,\frac{Y_j^2}{2}) =
Y_i Y_j Y_k -
a_{ik}\, Y_j \,\frac{Y_k^2}{2}- \\[.1in]
\qquad - a_{jk}\, Y_i\,\frac{Y_k^2}{2} + a_{ik}a_{jk} \,\frac{Y_k^3}{3}
-a_{ij}\,\frac{Y_j^2}{2} Y_k + a_{ij} a_{jk} Y_j \,\frac{Y_k^2}{2}
- a_{ij} a_{jk}^2 \frac{1}{2}\, \frac{Y_k^3}{3}.
\end{array}
$$

Let us fix $w\in W$ together with its reduced decomposition
$w=s_{i_1}\cdots s_{i_l}$.
Applying Lemma~\ref{lem:I_i(jjj)} repeatedly 
for the calculation of $\D_w = I_{i_l}\cdots I_{i_1}(1)$, 
and transferring the  sequences of integers $(k_1,\dots,k_n)$,
$n=1,2,\dots,l-1$, to the columns of a triangular array
$(k_{pq})$, we deduce the following result.

\begin{corollary} 
{\rm \cite{Du}}
For a reduced decomposition $w=s_{i_1}\cdots s_{i_l}\in W$, 
we have
$$
\D_w(y) = \sum_{(k_{pq})}\prod_{1\leq p<q\leq l}
\frac{(-a_{i_p i_q})^{k_{pq}}}{k_{pq}!}\,
\prod_{s=1}^l \frac{K_{*s}!\, Y_p^{K_{*s}+1 - K_{s*}}}{(K_{*s}+1 - K_{s*})!}\,,
$$
where the sum is over collections of nonnegative integers 
$(k_{pq})_{1\leq p<q\leq l}$ such that $K_{*s} + 1 \geq K_{s*}$,
for $s=1,\dots,l$; 
and $K_{*s} = \sum_{p<s} k_{ps}$ and $K_{s*} = \sum_{q>s} k_{sq}$.
\label{cor:duan}
\end{corollary}

This result is equivalent to Duan's recent result~\cite{Du}
about degrees $\deg(X_w)=\ell(w)!\,\D_w(\rho)$ of Schubert varieties.
Note that the approach and notations of~\cite{Du} 
are quite different from ours.

\section{K-theory and Demazure modules}
\label{sec:K-th-Demazure}

In this section, we recall a few facts about the K-theory for $G/B$.
\smallskip

Denote by $K(G/B)=K(G/B,\Q)$ the {\it Grothendieck ring\/} 
of coherent sheaves on $G/B$ with rational coefficients.
Let $\Q[\Lambda]$ be the group algebra of the weight lattice $\Lambda$.
It has a linear basis of formal exponents $\{e^\lambda\mid \lambda\in\Lambda\}$
with multiplication $e^\lambda\cdot e^\mu = e^{\lambda+\mu}$, i.e., 
$\Q[\Lambda]$ is the algebra of Laurent polynomials in the variables
$e^{\omega_1},\cdots,e^{\omega_r}$.  The action of the Weyl group on $\Lambda$
extends to a $W$-action on $\Q[\Lambda]$.
Let $\epsilon:\Q[\Lambda]\to \Q$ be the linear map such that 
$\epsilon(e^\lambda)=1$, for any $\lambda\in\Lambda$,
i.e., $\epsilon(f)$ is the sum of coefficients of exponents in $f$.
Then the Grothendieck ring $K(G/B)$ is canonically isomorphic\footnote
{The isomorphism is given by sending the K-theoretic
class $[\L_\lambda]_K\in K(G/B)$ of the line bundle $\L_\lambda$
to the coset $e^\lambda \pmod {\J_W}$, for any $\lambda\in\Lambda$.}
to the quotient ring:
$$
K(G/B)\simeq \Q[\Lambda]/\J_W,
$$
where $\J_W=\left<f\in \Q[\Lambda]^W\mid \epsilon(f)=0\right>$ is the ideal
generated by $W$-invariant elements $f\in\Q[\Lambda]$ with $\epsilon(f)=0$.
Let us identify the Grothendieck ring $K(G/B)$ with the quotient 
$\Q[\Lambda]/\J_W$ via this isomorphism.
Since $\epsilon$ annihilates the ideal $\J_W$, it induces
the map $\epsilon:K(G/B)\to\Q$, which we denote by the same letter.

The {\it Demazure operators\/}  
$T_i:\Q[\Lambda]\to \Q[\Lambda]$, $i=1,\dots,r$, are given by
\begin{equation}
T_i:f \mapsto  \frac{f - e^{-\alpha_i}\,s_i(f)}{1-e^{-\alpha_i}}.
\label{eq:demazure-operator}
\end{equation}
The Demazure operators satisfy the Coxeter relations 
$(T_i\, T_j)^{m_{ij}} = 1$ and $(T_i)^2 = T_i$.
For a reduced decomposition $w=s_{i_1}\cdots s_{i_l}\in W$,
define $T_w=T_{i_1}\cdots T_{i_l}$.  The operator $T_w$ depends only on 
$w\in W$ and does not depend on a choice of reduced decomposition.
The operators $T_i$ commute with operators of multiplication by 
$W$-invariant elements.  Thus the $T_i$ preserve the ideal $\J_W$ and induce 
operators acting on the Grothendieck ring $K(G/B)$, which we will 
denote by same symbols $T_i$.

The {\it Grothendieck classes\/} $\gamma_w\in K(G/B)$, $w\in W$, 
can be constructed, as follows.
$$
\begin{array}{l}
\displaystyle
\gamma_{w_\circ}=|W|^{-1}\, \prod_{\alpha\in\Phi^+} (1-e^{-\alpha}) 
\pmod {\J_W};\\[.15in]
\displaystyle
\gamma_w=T_{w^{-1}w_\circ}(\gamma_{w_\circ}),
\quad\textrm{for any }w\in W.
\end{array}
$$
According to Demazure~\cite{Dem}, the classes $\gamma_w$ 
are the K-theoretic classes $[\O_X]_K$ of the structure sheaves
of Schubert varieties $X=X_{w_\circ w}$.
In particular, $\gamma_{\id}=[\O_{G/B}]_K=1$.
The classes $\gamma_w$, $w\in W$, form a linear basis of $K(G/B)$.

Moreover, we have (see~\cite{Dem})
\begin{equation}
T_i(\gamma_{w}) = 
\left\{
\begin{array}{cl}
\gamma_{ws_i} &\textrm{if } \ell(ws_i)=\ell(w)-1,\\
\gamma_w &\textrm{if } \ell(ws_i)=\ell(w)+1.
\end{array}
\right.
\label{eq:T(gamma)}
\end{equation}




The {\it Chern character\/} is the ring isomorphism 
$\chi:K(G/B)\to H^*(G/B)$ 
induced by the map $\chi:e^\lambda \mapsto [e^\lambda]$,
for $\lambda\in\Lambda$, 
where $[e^\lambda]:=
1+\bar \lambda+\bar \lambda^2/2! + \cdots \in H^*(G/B)$
and $\bar\lambda=c_1(\L_\lambda)$, as before.
The isomorphism $\chi$ 
relates the Grothendieck classes $\gamma_w$ with the Schubert classes 
$\sigma_w$ by a triangular transformation:
\begin{equation}
\chi:\gamma_w\mapsto \sigma_w + \textrm{higher degree terms}.
\label{eq:Chern-char}
\end{equation}

For a dominant weight $\lambda\in\Lambda^+$,
let $V_\lambda$ denote the finite dimensional irreducible representation of 
the Lie group $G$ with highest weight $\lambda$. 
For $\lambda\in\Lambda^+$ and $w\in W$, the 
{\it Demazure module\/} $V_{\lambda,w}$ is the $B$-module 
that is dual to the space of global sections of the line bundle 
$\L_\lambda$ on the Schubert variety $X_w$:
$$
V_{\lambda,w} = H^0(X_w,\L_\lambda)^*.  
$$
For the longest Weyl group element $w=w_\circ$,
the space $V_{\lambda,w_\circ} = H^0(G/B,\L_\lambda)^*$ has the structure of
a $G$-module.  The classical {\it Borel-Weil theorem\/} says that 
$V_{\lambda,w_\circ}$ is isomorphic to the irreducible $G$-module $V_\lambda$.
Formal characters of Demazure modules are given by 
$\ch(V_{\lambda,w})=\sum_{\mu\in\Lambda} 
m_{\lambda,w}(\mu)\,e^\mu\in\Z[\Lambda]$,
where $m_{\lambda,w}(\mu)$ is the multiplicity of weight $\mu$ in 
$V_{\lambda,w}$.
They generalize characters of irreducible representations
$\ch(V_\lambda)=ch(V_{\lambda,w_\circ})$.
{\it Demazure's character formula}~\cite{Dem} says that 
the character $ch(V_{\lambda,w})$ is given by
\begin{equation}
\ch(V_{\lambda,w}) = T_w(e^\lambda).
\label{eq:demazure-character-formula}
\end{equation}

\section{Asymptotic expression for degree}
\label{sec:asymtotic}

\begin{proposition}
For any $w\in W$, the dimension of the Demazure module $V_{\lambda,w}$
is a polynomial in $\lambda$ of degree $\ell(w)$.
The polynomial $\D_w$ is the leading homogeneous component
of the polynomial $\dim V_{\lambda,w}\in \Q[V]$.
In other words, the value $\D_w(\lambda)$ equals
$$
\D_w(\lambda) = \lim_{k\to\infty} \frac{\dim V_{k\lambda,w}}{k^{\ell(w)}},
$$
for any $\lambda\in\Lambda^+$.
\label{prop:demazure-limit}
\end{proposition}

Proposition~\ref{prop:demazure-limit} together with Weyl's dimension
formula implies the following statement, which was derived by
Stembridge using Standard Monomial Theory.

\begin{corollary}  {\rm \cite[Theorem~1.1]{Ste}} \
For the longest Weyl group element $w=w_\circ$, we have
$$
\D_{w_\circ}(y) = \prod_{\alpha\in\Phi^+} 
\frac {(y,\alpha^\vee)}{(\rho,\alpha^\vee)}.
$$
\label{cor:D-w0}
\end{corollary}

\begin{proof}  
Weyl's formula says that
the dimension of $V_{\lambda,w_\circ}=V_\lambda$ is
$$
\dim V_\lambda = \prod_{\alpha\in\Phi^+} 
\frac {(\lambda+\rho,\alpha^\vee)}{(\rho,\alpha^\vee)}.
$$
Taking the leading homogeneous component of this polynomial in $\lambda$, 
we prove the claim for $y=\lambda\in\Lambda^+$, and thus, for any $y\in V$.
\end{proof}

In order to prove Proposition~\ref{prop:demazure-limit} we need
the following lemma.

\begin{lemma}  The map $\epsilon:K(G/B)\to\Q$
is given by $\epsilon(\gamma_w) = \delta_{w,\id}$,
for any $w\in W$.
\label{lem:epsilon}
\end{lemma}

\begin{proof}
It follows directly from the definitions that
the Chern character $\chi$ translates $\epsilon$ to 
the map $\epsilon\cdot\chi^{-1}:H^*(G/B)\to\Q$ 
given by $\epsilon\cdot \chi^{-1}:\bar f\mapsto f(0)$, for a polynomial
representative $f\in\Q[\h]$ of $\bar f$.
Thus $\epsilon\cdot\chi^{-1}(\sigma_w)=\delta_{w,\id}$. 
Indeed, $\sigma_{\id}=1$ and all other Schubert classes $\sigma_w$ have zero 
constant term, for $w\ne\id$.
Triangularity~(\ref{eq:Chern-char}) of the Chern character
implies the needed statement.
\end{proof}

\begin{proof}[Proof of Proposition~\ref{prop:demazure-limit}]
The preimage of identity~(\ref{eq:e-lambda-D-uw}), for $u=\id$,
under the Chern character $\chi$ is the following 
expression in $K(G/B)$: 
$$
e^\lambda=\sum_{w\in W} \D_w(\lambda)\,\chi^{-1}(\sigma_w)
=\sum_{w\in W} \hat \D_w(\lambda)\, \gamma_w,
$$
Triangularity~(\ref{eq:Chern-char}) implies that
$\chi^{-1}(\sigma_w)=\gamma_w+ \sum_{\ell(u)>\ell(w)} 
c_{w,u}\, \gamma_u$ and 
$\hat \D_w = \D_w + \sum_{\ell(u)<\ell(w)} c_{u,w} \,\D_u$,
for some coefficients $c_{w,u}\in\Q$.
Thus the homogeneous polynomial $\D_w$ is the leading homogeneous component of 
the polynomial $\hat\D_w$.
Applying  the map $\epsilon \cdot T_w$ to both sides of the previous expression
and using Lemma~\ref{lem:epsilon}, we obtain
$$
\epsilon (T_w (e^\lambda))= \sum_{u\leq w} \hat\D_u(\lambda).
$$
Indeed, according to~(\ref{eq:T(gamma)}),
the coefficient of $\gamma_{\id}$ in $T_w(\gamma_u)$ 
is equal to 1 if $u\leq w$, and 0 otherwise.  Thus 
$\epsilon (T_w (e^\lambda))$ is a  polynomial in $\lambda$
of degree $\ell(w)$
and its leading homogeneous component is again $\D_w$.  
But, Demazure's character formula says that 
$T_w (e^\lambda)$ is the character of $V_{\lambda,w}$ and
$\epsilon (T_w (e^\lambda))=\dim V_{\lambda,w}$.
\end{proof}

Lakshmibai reported the following simple geometric proof of 
Proposition~\ref{prop:demazure-limit}.
Assume that $\lambda$ is a dominant regular weight.
We have $V_{w,k\lambda}^*=H^0(X_w,\L_{k\lambda}) = R_k$,
where $R_k$ is the $k$-th graded component of the coordinate ring $R$
of the image of $X_w$ in $\mathbb{P}(V_\lambda)$.  The Hilbert polynomial
of the coordinate ring has the form $\mathrm{Hilb}_R(k) = \dim R_k =
A\, k^l/l! + (\textrm{lower degree terms})$, where $l=\dim_\C X_w = \ell(w)$, 
and $A=\deg_\lambda(X_w)$ is the degree of $X_w$ in $\mathbb{P}(V_\lambda)$.
Thus $\lim_{k\to\infty} \dim V_{k\lambda,w}/k^{\ell(w)} = A/l!=
\deg_\lambda(X_w)/ \ell(w)! = \D_w(\lambda)$.

\section{Permanent of the matrix of Cartan integers}
\label{sec:permanent}

Let us give a curious consequence of Theorem~\ref{th:D=S}.

\begin{corollary}
Let $A=(a_{\alpha,\beta})$ be the $N\times N$-matrix, $N=|\Phi^+|$,
formed by the Cartan integers $a_{\alpha,\beta}=(\alpha,\beta^\vee)$,
for $\alpha,\beta\in\Phi^+$.
Then the permanent of the matrix $A$ equals
$$
\mathrm{per} (A) = |W|\cdot\prod_{\alpha\in\Phi^+} (\rho,\alpha^\vee).
$$
\label{cor:permanent}
\end{corollary}

The matrix $A$ should not be confused with the Cartan matrix.  The latter
is a certain $r\times r$-submatrix of $A$.

\begin{proof}
According to Theorem~\ref{th:D=S} and Corollary~\ref{cor:D-w0}, 
we have  
$$
1=\D_{\id}=\sigma_{w_\circ}(\ddy)\cdot \D_{w_\circ}(y)=
\left(\frac{1}{|W|}\,\prod_{\alpha\in\Phi^+} d_\alpha\right)\cdot
\left(\prod_{\beta\in\Phi^+} \frac{(y,\beta^\vee)}{(\rho,\beta^\vee)}
\right),
$$
where $d_\alpha$ is the operator of differentiation with respect 
to a root $\alpha$ given by~(\ref{eq:diff-mu}).
Using the product rule for differentiation and 
the fact that $d_\alpha\cdot(y,\beta^\vee)=
(\alpha,\beta^\vee)$, we derive the claim.
\end{proof}

For type $A_{n-1}$, we obtain the following result. 

\begin{corollary}
Let $B=(b_{ij,k})$ be the $\binom n2\times n$-matrix with rows labeled
by pairs $1\leq i<j\leq n$ and columns labeled by $k=1,\dots,n$ such 
that $b_{ij,k} = \delta_{i,j}-\delta_{j,k}$.
Then
$$
\mathrm{per}(B\cdot B^T) = 1!\,2!\cdots n!.
$$
\label{cor:permanent-A}
\end{corollary}

\begin{proof}
For type $A_{n-1}$, the matrix $A$ in Corollary~\ref{cor:permanent}
equals $B\cdot B^T$.
\end{proof}

This claim can be also derived from the Cauchy-Binet formula
for permanents.

For example, for type $A_3$, we have
$$
\mathrm{per}
\left(
\left[
\begin{matrix}
1&-1&0&0\\
1&0&-1&0\\
1&0&0&-1\\
0&1&-1&0\\
0&1&0&-1\\
0&0&1&-1
\end{matrix}
\right]
\cdot
\left[
\begin{matrix}
1 & 1& 1& 0 & 0 & 0 \\
-1 & 0& 0& 1 & 1 & 0 \\
0 & -1& 0& -1 & 0 & 1 \\
0 & 0& -1& 0 & -1 & -1 
\end{matrix}
\right]
\right)
= 1! \, 2! \, 3!\, 4!.
$$

Note that the rank of the $\binom n2\times \binom n2$-matrix $B\cdot B^T$
is at most $n-1$.  Thus the determinant of this matrix is zero, 
for $n\geq 3$.
It would be interesting to find a combinatorial proof of 
Corollary~\ref{cor:permanent-A}.

\section{Schubert polynomials}
\label{sec:A-Schubert-poly}

In the rest of the paper we will be mainly concerned
with the case $G=SL_n$.  

The root system $\Phi$ associated to $SL_n$ is of the type $A_{n-1}$.
In this case, the spaces $V$
can be presented as 
$V=\Q^n/(1,\dots,1)\Q$.
Then $\Phi=\{\varepsilon_i-\varepsilon_j\in V\mid 1\leq i\ne j\leq n\}$,
where the $\varepsilon_i$ are images of the coordinate vectors in $\Q^n$.
The Weyl group is the symmetric group $W=S_n$ of order $n$ that acts on
$V$ 
by permuting the coordinates in $\Q^n$.
The Coxeter generators are the adjacent transpositions $s_i=(i,i+1)$.
The length $\ell(w)$ of a permutation $w\in S_n$ is the number of
inversions in $w$.
The longest permutation in $S_n$ is $w_\circ = n,n-1,\cdots ,2,1$.


The quotient $SL_n/B$ is the classical complex flag variety.
Its cohomology ring $H^*(SL_n/B)$ over $\Q$ is canonically
identified with the quotient
$$
H^*(SL_n/B) = \Q[x_1,\dots,x_n]/\I_n,
$$
where $\I_n=\<e_1,\dots,e_n\>$ is the ideal generated by the 
elementary symmetric polynomials $e_i$ in the variables $x_1,\dots,x_n$.
The divided difference operators $A_i$ act on the polynomial
ring $\Q[x_1,\dots,x_n]$ by
$$
A_i:f(x_1,\dots,x_n)\mapsto 
\frac{f(x_1,\dots,x_n)-f(x_1,\dots,x_{i-1},x_{i+1},x_i,x_{i+1},\dots,x_n)}{x_i-x_{i+1}}.
$$
For a reduced decomposition $w=s_{i_1}\cdots s_{i_l}$,
let $A_w=A_{i_1}\cdots A_{i_l}$.

Lascoux and Sch\"utzenberger~\cite{LS} defined the 
{\it Schubert polynomials\/} $\S_w$, for $w\in S_n$, by
$$
\S_{w_\circ} = x_1^{n-1} x_2^{n-2}\cdots x_{n-1}
\quad\textrm{and}\quad
\S_w = A_{w^{-1}w_\circ}(\S_{w_\circ}).
$$
Then the cosets of Schubert polynomials $\S_w$ modulo 
the ideal $\I_n$ are the Schubert classes $\sigma_w=\bar\S_w$ in 
$H^*(SL_n/B)$.

This particular choice of polynomial representatives
for the Schubert classes has the following {\it stability property.}
The symmetric group $S_n$ is naturally embedded into $S_{n+1}$ as the 
set of order $n+1$ permutations that fix the element $n+1$.
Then the Schubert polynomials remain the same under this embedding.

Let $S_\infty$ be the injective limit of symmetric groups
$S_1\hookrightarrow S_2\hookrightarrow S_3\hookrightarrow\cdots$.
In other words, $S_\infty$ is the group of infinite permutations 
$w:\Z_{>0}\to\Z_{>0}$ such that $w(i)=i$ for almost all $i$'s.
We think of $S_n$ as the subgroup of infinite permutations $w\in S_\infty$
that fix all $i>n$.
Let $\Q[x_1,x_2,\dots]$ be the polynomial ring in infinitely many variables
$x_1,x_2,\dots$.  The stability of the Schubert polynomials under the 
embedding $S_n\hookrightarrow S_{n+1}$ implies that the Schubert polynomials
$\S_w\in\Q[x_1,x_2,\dots]$ are consistently defined for any $w\in S_\infty$.
Moreover,
$\{\S_w\}_{w\in S_\infty}$ is a basis of the polynomial ring 
$\Q[x_1,x_2,\dots]$.

\section{Degree polynomials for type $A$}
\label{sec:degree-type-A}

Let us summarize properties of the polynomials $\D_{u,w}$
for type $A_{n-1}$.

Let $y_1,\dots,y_n$ be independent variables.
Let us assign to each edge $w\lessdot w s_{ij}$ in the Hasse diagram 
of the Bruhat order on $S_n$ the weight $m(w, w s_{ij})=
y_i-y_j$.  For a saturated chain
$C=(u_0\lessdot u_1\lessdot u_2\lessdot\cdots\lessdot u_l)$
in the Bruhat order, we define its weight as
$m_C(y)=m(u_0,u_1)\,m(u_1,u_2)\cdots m(u_{l-1},u_l)$.

For $u, w\in S_n$ such that $u\leq w$, the polynomial 
$\D_{u,w}\in\Q[y_1,\dots,y_n]$
is defined as the sum
$$
\D_{u,w} = \frac{1}{\ell(w)!}\,\sum_{C} m_C(y)
$$
over all saturated chains $C=(u_0\lessdot u_1 \lessdot \cdots \lessdot u_l)$
from $u_0=u$ to $u_l=w$ in the Bruhat order.
Also $\D_{w}:=\D_{\id, w}$.

The subspace $\H_n$ of $S_n$-harmonic polynomials 
in $\Q[y_1,\dots,y_n]$ is given by 
$$
\H_n=\{g\in\Q[y_1,\dots,y_n]\mid
f(\ddy_1,\dots,\ddy_n)\cdot g(y_1,\dots,y_n)=0\textrm{ for any }f\in\I_n\}.
$$

\begin{corollary}
{\rm (1)}
The polynomials $\D_w$, $w\in S_n$,
form a basis of $\H_n$.

{\rm (2)}
The polynomials $\D_{u,w}$, $u,w\in S_n$, can be expressed as
$$
\begin{array}{l}
\displaystyle
\D_{w_\circ} = \frac 1 {1!\,2!\cdots (n-1)!}
\prod_{1\leq i<j\leq n}(y_i-y_j)
=\det\left(\left(y_i^{(n-j)}\right)_{i,j=1}^n\right),
\\[.2in]
\D_{u,w}=\S_u(\ddy_1,\dots,\ddy_n)\,\S_{w_\circ w}(\ddy_1,\dots,\ddy_n)\cdot\D_{w_\circ},
\end{array}
$$
where $a^{(b)} = \frac{a^b}{b!}$.

{\rm (3)}
The polynomials $\D_w$, $w\in S_n$, can be also expressed as
$$
\D_w=I_{w^{-1}}(1),
$$
where $I_{w}=I_{i_1}\cdots I_{i_l}$,
for a reduced decomposition $w=s_{i_1}\cdots s_{i_l}$,
and the operators $I_1,\dots,I_{n-1}$ on 
$\Q[y_1,\dots,y_n]$ are given by
$$
I_i:g(y_1,\dots,y_n)\mapsto\int_0^{y_i-y_{i+1}}
g(y_1,\dots,y_{i-1},y_i-t,y_{i+1}+t,
y_{i+2},\dots,y_n)\,dt.
$$
\label{cor:D-S}
\end{corollary}

The following symmetries are immediate from the definition 
of the polynomials~$\D_w$.

\begin{lemma} 
{\rm (1)}
For any $w\in S_n$, we have
$$
\D_w(y_1,\dots,y_n)=\D_{w_\circ w w_\circ}(-y_n,\dots,-y_1).
$$

{\rm(2)} Also $\D_w(y_1+c,\dots,y_n+c) = \D_w(y_1,\dots,y_n)$,
for any constant $c$.
\label{lem:involution-omega}
\end{lemma}

The spaces $\H_n$ of $S_n$-harmonic polynomials are embedded
in the polynomial ring $\Q[y_1,y_2,\dots]$ in infinitely many variables: 
$\H_1\subset \H_2\subset \H_3\subset \cdots \subset \Q[y_1,y_2,\dots]$.
Moreover, the union of all $\H_n$'s is exactly this polynomial ring.
It is clear from the definition that the polynomials $\D_w$ 
are stable under the embedding $S_n\hookrightarrow S_{n+1}$.
Thus the polynomials $\D_w\in\Q[y_1,y_2,\dots]$
are consistently defined for any $w\in S_\infty$.

\begin{corollary}
{\rm (1)}
The set of polynomials $\D_w$, $w\in S_\infty$, forms a linear basis 
of the polynomial ring $\Q[y_1,y_2,\dots]$. 

{\rm (2)}
The basis $\{\S_w\}_{w\in S_\infty}$ of Schubert polynomials in 
$\Q[x_1,x_2,\dots]$ is D-dual\/%
\footnote{Note that $D$-pairing between polynomials in $n$ variables
is stable under the embedding $\Q[x_1,\dots,x_n]\subset 
\Q[x_1,\dots,x_{n+1}]$.  Thus $D$-pairing is consistently defined
for polynomials in infinitely many variables.}
to the basis 
$\{\D_w\}_{w\in S_\infty}$ in $\Q[y_1,y_2,\dots]$,
i.e.,  $(\S_u,\D_w)_D=\delta_{u,w}$, for any $u,w\in S_\infty$.
\end{corollary}

\begin{proof}   Let $u,v\in S_\infty$.  Then, for sufficiently
large $n$, we have  $u,v\in S_n$.  Now the identity 
$(\S_u,\D_w)_D=\delta_{u,w}$ follows from Corollary~\ref{cor:D-dual-sigma}.
\end{proof}

\section{Flagged Schur polynomials}
\label{sec:flagged-Schur}

Let $\mu=(\mu_1,\dots,\mu_n)$, $\mu_1\geq \cdots\geq \mu_n \geq 0$, 
be a partition, $\beta=(\beta_1,\dots,\beta_m)$ be a nonnegative
integer sequence, and  
$a=(a_1\leq \dots\leq a_n)$ and $b=(b_1\leq \dots\leq b_n)$ 
be two weakly increasing positive integer sequences.
A {\it flagged semistandard Young tableau\/} 
of shape $\mu$, weight $\beta$, with flags $a$ and $b$ 
is an array of positive integers $T=(t_{ij})$, $i=1,\dots,n$, 
$j=1,\dots,\mu_i$,
such that
\begin{enumerate}
\item entries strictly increase in the columns: $t_{1j}<t_{2j}<t_{3j} <\cdots$;
\item entries weakly increase in the rows: $t_{i1}\leq t_{i2}\leq t_{i3} <\cdots$;
\item $\beta_k=\#\{(i,j)\mid t_{ij} = k\}$ is the number of entries $k$ in $T$,
for $k=1,\dots,m$;
\item for all entries in the $i$-th row, we have $a_i\leq t_{ij}\leq b_i$.
\end{enumerate}
The {\it flagged Schur polynomial\/} $s_\mu^{a,b}=s_\mu^{a,b}(x)
\in\Q[x_1,x_2,\dots]$ is defined as the sum
$$
s_\mu^{a,b}(x) = \sum x^T
$$
over all flagged semistandard Young tableaux $T$ of shape $\mu$ 
with flags $a$ and $b$ and arbitrary weight, 
where $x^T := x_1^{\beta_1} \cdots x_m^{\beta_m}$
and $\beta=(\beta_1,\dots,\beta_m)$ is the weight of $T$.

Note that $s_\mu^{(1,\dots,1),(n,\dots,n)}$ is the usual Schur
polynomial $s_\mu(x_1,\dots,x_n)$.
Flagged Schur polynomials were originally introduced by 
Lascoux and Sch\"utzenberger~\cite{LS}.

The polynomial $s_\mu^{a,b}(x)$ does not depend on the flag $a$ provided that
$a_i\leq i$, for $i=1,\dots,n$.  Indeed, entries in the $i$-th row of 
any semistandard Young tableaux (of a standard shape) are greater
than or equal to $i$.  Thus the condition $a_i\leq t_{ij}$ is redundant.
Let 
$$
s_\mu^b(x) := s_\mu^{(1,\dots,1),b}(x)=s_\mu^{(1,\dots,n),b}(x).
$$

Flagged semistandard Young tableaux can be presented
as collections of $n$ non-crossing lattice paths on $\Z\times\Z$
that connect points $A_1,\dots,A_n$ with $B_1,\dots,B_n$,
where $A_i =(-i, a_i)$ and $B_i=(\mu_i-i,b_i)$.
Let us assign the weight $x_i$ to each edge $(i,j)\to(i,j+1)$ in a lattice 
path and weight $1$ to an edge $(i,j)\to(i+1,j)$.
Then the product of weights over all edges in the collection of lattice
paths corresponding to a flagged tableau $T$ equals $x^T$.
According to the method of Gessel and Viennot~\cite{GV} for 
counting non-crossing lattice paths, the flagged Schur polynomial 
$s_\mu^{a,b}(x)$ equals the determinant
\begin{equation}
s_\mu^{a,b}(x) = \det \left( h_{\mu_i-i+j}^{[a_j,b_i]}\right)_{i,j=1}^n,
\label{eq:flagged-schur=det}
\end{equation}
where, for $k\leq l$,
$$
h_m^{[k,l]} = h_m(x_k,x_{k+1},\dots,x_l)=
\sum_{k\leq i_1\leq \cdots \leq i_m\leq l} x_{i_1}\cdots x_{i_m}
$$
is the complete homogeneous symmetric polynomial of degree $m$ in the variables 
$x_k,\dots,x_l$;  and $h_m^{[k,l]}=0$, for $k>l$.
Another proof of this result was given by Wachs~\cite{Wac}.

For permutations $w=w_1\cdots w_n$ in $S_n$ and
$\sigma=\sigma_1\cdots \sigma_r$ in $S_r$, let us say 
that $w$ is {\it $\sigma$-avoiding\/} if there
is no subset $I=\{i_1<\cdots < i_r\}\subseteq \{1,\dots,n\}$
such that the numbers $w_{i_1},\dots,w_{i_r}$ have the same relative order 
as the numbers $\sigma_1,\dots,\sigma_r$.
Let $S_n^\sigma\subseteq S_n$ be the set of $\sigma$-avoiding permutations 
in $S_n$.
For example, a permutation $w=w_1\cdots w_n$ is {\it 312-avoiding\/} 
if there are no $i<j<k$ such that $w_i>w_k>w_j$.   
It is well-known that, for any permutation $\sigma\in S_3$ of size 3, 
the number of $\sigma$-avoiding permutations in $S_n$ equals the 
Catalan number $\frac{1}{n+1} \binom{2n}{n}$.
A permutation $w$ is called {\it vexillary\/} if it is 
2143-avoiding.

Lascoux and Sch\"utzenberger~\cite{LS} stated that
Schubert polynomials for vexillary permutations are certain
flagged Schur polynomials.  This claim was clarified and proved by 
Wachs~\cite{Wac}.

For a permutation $w=w_1\cdots w_n$ is $S_n$, the inversion sets
$\Inv_i(w)$, $i=1,\dots,n$, are defined as 
$$
\Inv_i(w) = \{j\mid i<j\leq n \textrm{ and } w_i>w_j\}.
$$
The {\it code\/} of permutation $w$ is the sequence
$\code(w)=(c_1,\dots,c_n)$ given by
$$
c_i=c_i(w)=|\Inv_i(w)|=\#\{j\mid j>i,\ w_j<w_i\}
\textrm{ for } i=1,\dots,n.
$$
The map $w\mapsto \code(w)$ is a bijection
between the set of permutations $S_n$ and the set of vectors 
$\{(c_1,\dots,c_n)\in\Z^n\mid
0\leq c_i\leq n-i,\textrm{ for }i=1,\dots,n\}$.

The {\it shape\/} of permutation $w\in S_n$ is the partition 
$\mu=(\mu_1\geq \dots\geq \mu_m)$ 
given by nonzero components $c_i$ of its code 
arranged in decreasing order.
The {\it flag\/} of permutation $w\in S_n$ is the sequence
$b=(b_1\leq \cdots \leq b_m)$ given by
the numbers $\min \Inv_i(w) -1$, for non-empty $\Inv_i(w)$, arranged 
in increasing order.

\begin{proposition} {\rm \cite{Wac}, cf.~\cite{LS}} \
Assume that $w\in S_n^{2143}$ is a vexillary permutation. 
Let $\mu$ be its shape and $b$ be its flag.
Then the Schubert polynomial $\S_w(x)$ is the following flagged 
Schur polynomial:
$\S_w(x) = s_\mu^b(x)$.
\label{prop:Schubert=Schur}
\end{proposition}

We remark that not every flagged Schur polynomial is a Schubert
polynomial.

Let $\CC_n$ be the set of partitions $\mu=(\mu_1,\dots,\mu_n)$,
$\mu_1\geq \cdots\geq \mu_n\geq 0$, such that $\mu_i\leq n-i$,
for $i=1,\dots,n$, i.e., $\CC_n$ is the set of partitions 
whose Young diagrams fit inside the staircase shape
$(n-1,n-2,\dots,0)$.  These partitions are in an obvious correspondence
with Catalan paths.  Thus $|\CC_n|=\frac{1}{n+1}\binom{2n}{n}$ is the Catalan number.

A permutation $w$ is called {\it dominant\/} if 
$\code(w)=(c_1,\dots,c_n)$ is a partition, i.e., $c_1\geq \cdots 
\geq c_n$.  The next claim is essentially well known;
see, e.g.,~\cite{Man}.

\begin{proposition}
A permutation $w=w_1\cdots w_n\in S_n$ is dominant if and only if it
is 132-avoiding.

The map $w\mapsto \code(w)$ is a bijection between the set $S^{132}_n$ of
dominant permutations and the set $\CC_n$.
We have $w_i>w_{i+1}$ if and only if $c_i>c_{i+1}$,
and $w_i<w_{i+1}$ if and only if $c_i=c_{i+1}$.

For $w\in S_n^{132}$, we have
$\Inv_i(w) = \{k\mid w_k < \min\{w_1,\dots,w_i\}\}$ and
$c_i(w) = \min\{w_1,\dots,w_i\} -1$.

The inverse map $c\mapsto w(c)$ from $\CC_n$ to $S_n^{132}$
is given recursively by $w_1=c_1+1$ and 
$w_i=\min\{j>c_i\mid j\ne w_1,\dots,w_{i-1}\}$, for $i=2,\dots,n$.
In particular, if $c_{i}<c_{i-1}$ then $w_i=c_i+1$.
\label{prop:132-avoiding}
\end{proposition}

\begin{proof}
Let us assume that $w$ is 132-avoiding and show that $\code(w)$
is weakly decreasing.
Indeed, if $w_i>w_{i+1}$ then $c_i>c_{i+1}$.
If $w_i<w_{i+1}$ then there is no $j>i+1$ such that $w_i<w_j<w_{i+1}$,
because $w$ is 132-avoiding.  Thus $c_i=c_{i+1}$ in this case.

On the other hand, assume that $w\in S_n$ is not a 132-avoiding permutation.
Say that $(i,j,k)$ is a 132-triple of indices if $i<j<k$ and $w_i<w_k<w_j$.
Let us find a 132-triple $(i,j,k)$ such that the difference $j-i$ is 
as small as possible.  We argue that $j=i+1$.  Otherwise, pick any $l$ such
that $i<l<j$.  If $w_l<w_k$ then $(l,j,k)$ is a 132-triple, and 
if $w_l>w_k$ then $(i,l,k)$ is a 132-triple.  Both these triples have
a smaller difference.  This shows that we can always find a 132-triple 
of the form $(i,i+1,k)$.  Then $c_i(w)<c_{i+1}(w)$.
Thus  $\code(w)$ is not weakly decreasing.  This proves that
$w\mapsto \code(w)$ is a bijection between $S_n^{132}$ and $\CC_n$.

Let $w\in S_n^{132}$.  Fix an index $i$ and find $1\leq j\leq i$ such 
that $w_j=\min\{w_1,\dots,w_i\}$.
Since $w$ is 132-avoiding, there is no $k>i$ such that $w_i>w_k>w_j$.
Thus the conditions $k>i$, $w_k<w_i$ imply that $w_k<w_j$.
On the other hand, if $w_k<w_j$ for some $k\in\{1,\dots,n\}$ then $k>i$
because of our choice of $j$.  This shows that the $i$-th inversion set
of the permutation $w$ is
$\Inv_i(w) = \{k\mid w_k < \min\{w_1,\dots,w_i\}\}$.
Thus $c_i(w) = |\Inv_i(w)|=\min\{w_1,\dots,w_i\}-1$.

Let $w\in S_n^{132}$ and $\code(w)=(c_1,\dots,c_n)$.  We have $w_1=c_1+1$.  
Let us derive the identity $w_{i}=\min\{j>c_{i}\mid j\ne w_1,\dots,w_{i-1}\}$,
for $i=2,\dots,n$.  Indeed, if $c_i<c_{i-1}$ then $w_i=c_i+1$, as needed.
Otherwise, if $c_i=c_{i-1}$, then $w_i> w_{i-1}$.
Let $k$ be the index such that $w_{k}=
\min\{j>c_{i}\mid j\ne w_1,\dots,w_{i-1}\}$.  If $k\ne i$ then
$k>i$ and $w_k<w_i$.  Thus $w_{i-1}<w_k<w_i$.  This is impossible because
we assumed that $w$ is $132$-avoiding.
\end{proof}

The following claim is also well known; see, e.g.,~\cite{Man}.

\begin{proposition}
For a dominant permutation $w\in S_n^{132}$, the Schubert polynomial
is given by the monomial $\S_w(x) = x_1^{c_1(w)}\cdots x_n^{c_n(w)}$. 
\label{cor:S=monimial}
\end{proposition}

This claim follows from Proposition~\ref{prop:Schubert=Schur},
because the set of dominant permutations is a subset of vexillary 
permutations.  

\begin{proof}
Let $\mu=\code(w) = (k_1^{m_1},k_2^{m_2},\dots)$,  $k_1>k_2>\cdots$,
be the shape of $w$.  According to Proposition~\ref{prop:132-avoiding},
the flag of $w$ is $b=(m_1^{m_1},(m_1+m_2)^{m_2},\dots)$.
For this shape and flag, there exists only one flagged semistandard Young 
tableau $T=(t_{ij})$,
which is given by $t_{ij}=i$.
Thus $\S_w(x) = s_\mu^b(x)=x_1^{\mu_1}\cdots x_n^{\mu_m}$.
\end{proof}

A permutation $w$ is 3412-avoiding if and only if
$w_\circ w$ is vexillary.  Also a permutation $w$ is 312-avoiding 
if and only if $w_\circ w$ is 132-avoiding.
The next claim follows from Theorem~\ref{th:D=S}, 
Proposition~\ref{prop:Schubert=Schur}, and Corollary~\ref{cor:S=monimial}.

\begin{theorem}
Let $w\in S_n^{3412}$ be a 3412-avoiding permutation.  Let
$\mu$ and $b$ be the shape and flag of the vexillary
permutation $w_\circ w$.  Then 
$$
\D_w(y_1,\dots,y_n) = \frac{1}{1!\,2!\,\cdots (n-1)!}\,
s_\mu^b(\ddy_1,\dots,\ddy_n)\cdot \prod_{i<j}(y_i-y_j).
$$

In particular, for a 312-avoiding permutation $w\in S_n^{312}$ 
and $(c_1,\dots,c_n)=\code(w_\circ w)$, we have  
$$
\begin{array}{l}
\displaystyle
\D_w(y_1,\dots,y_n)= \frac{1}{1!\,2!\,\cdots (n-1)!}\,
\left(\prod_{k=1}^n(\ddy_k)^{c_k}\right) \cdot \prod_{i<j}
(y_i-y_j)\\[.1in]
\qquad \qquad \qquad 
\displaystyle
=\det\left( \left(y_i^{(n-c_i-j)}\right)_{i,j=1}^n\right),
\end{array}
$$
where $a^{(b)}=a^b/b!$, for $b\geq 0$, and $a^{(b)}=0$, for $b<0$.
\label{th:3412-avoiding}
\end{theorem}


Applying Lemma~\ref{lem:involution-omega}(1), we obtain the determinant 
expression for $\D_w$, for 231-avoiding permutations $w$, as well.

\begin{corollary} For a 231-avoiding permutation $w\in S_n^{231}$ 
and $(c_1,\dots,c_n)=\code(w w_\circ)$, we have  
$$
\D_w(y_1,\dots,y_n)= 
\det\left( \left((-y_{n-i+1})^{(n-c_i-j)}\right)_{i,j=1}^n\right).
$$
\end{corollary}

\section{Demazure characters for 312-avoiding permutations}
\label{sec:demazure-312}

In the previous section we gave a simple determinant formula for 
the polynomial $\D_w$, for a 312-avoiding permutation $w\in S_n^{312}$.
We remark that 312-avoiding permutations are exactly the {\it Kempf elements\/}
that were studied by Lakshmibai in~\cite{Lak}.
In this and the following sections, we give some additional nice
properties of 312-avoiding permutations.
In this section, we show how Weyl's character formula can be easily 
deduced from Demazure's character formula by induction on some 
sequence of 312-avoiding permutations that interpolates between 1
and $w_\circ$.

\smallskip

Let $z_1,\dots,z_n$ be independent variables, and let 
$T_i$, $i=1,\dots,n-1$, be the operator that acts on 
the polynomial ring $\Q[z_1,\dots,z_n]$ by
$$
T_i:f(z_1,\dots,z_n)\mapsto  \frac{z_i\,f(z_1,\dots,z_n)-
z_{i+1}\,f(z_1,\dots,z_{i-1},z_{i+1},z_i,z_{i+2},\dots,z_n)}
{z_i-z_{i+1}}.
$$
For $\lambda=(\lambda_1\geq \cdots\geq \lambda_n)$ and 
a reduced decomposition $w=s_{i_1}\cdots s_{i_l}\in S_n$, 
let
$$
\ch_{\lambda,w}(z_1,\dots,z_n) = T_{i_1}\cdots T_{i_l}(z_1^{\lambda_1}\cdots z_n^{\lambda_n}).
$$
The polynomials $\ch_{\lambda,w}$ do not depend on choice of 
reduced decomposition for $w$ because the $T_i$ satisfy the Coxeter 
relations.  Let us map the ring $\Q[z_1,\dots,z_n]$ to the 
group algebra $\Q[\Lambda]$ of the type $A_{n-1}$ weight lattice $\Lambda$
by $z_i\mapsto e^{\omega_i-\omega_{i-1}}$, for $i=1,\dots,n$,
where we assume that $\omega_0=\omega_{n}=0$.
Then the operators $T_i$ specialize to the 
Demazure operators~(\ref{eq:demazure-operator}) and the 
polynomials $\ch_{\lambda,w}$ map to the characters of Demazure modules
$\ch(V_{\lambda,w})$; cf.\ the Demazure character 
formula~(\ref{eq:demazure-character-formula}).
The polynomials $\ch_{\lambda,w}$ were studied by
Lascoux and Sch\"utzenberger~\cite{LS}, who called them
{\it essential polynomials\/}, and  
by Reiner and Shimozono~\cite{RS}, who called them
{\it key polynomials.}
To avoid confusion, we will call the polynomials $\ch_{\lambda,w}$ 
simply {\it Demazure characters.}

For a given partition $\lambda=(\lambda_1,\dots,\lambda_n)$,
the number of nonzero flagged Schur polynomials 
$s_\lambda^b(z_1,\dots,z_n)$ in $n$ variables equals 
the Catalan number $\frac{1}{n+1} \binom{2n}{n}$.
Indeed, such a polynomial is nonzero if and only if the flag 
$b=(b_1,\dots,b_n)$ satisfies $b_1\leq \dots\leq b_n\leq n$ and 
$b_i\geq i$, for $i=1,\dots,n$.
Let us denote by $\tilde\CC_n$ the set of such flags $b$.
The map $(b_1,\dots,b_n)\mapsto (c_1,\dots,c_n)$ given by
$c_i=n-b_i$, for $i=1,\dots,n$, is a bijection between
the sets $\tilde\CC_n$ and $\CC_n$.
The next theorem says that the flagged Schur polynomials 
$s_\lambda^b(z_1,\dots,z_n)$ are exactly 
the Demazure characters $\ch_{\lambda,w}$, for 312-avoiding
permutations $w\in S_n$.

Recall that the map $w\mapsto \code(w)$ is a bijection
between the sets $S_n^{132}$ and $\CC_n$
(see Proposition~\ref{prop:132-avoiding}).
Then the map $w\mapsto b(w)=(b_1,\dots,b_n)$ given by 
$b_i=n-c_i(w_\circ w)$, for $i=1,\dots,n$, 
is a bijection between the sets $S_n^{312}$ and $\tilde\CC_n$.
Note that $\ell(w)=b_1+\cdots+b_n-\binom {n+1}2$.
The inverse map $\tilde\CC_n\to S_n^{312}$
can be described recursively, as follows: 
$w_1=b_1$ and $w_i=\max\{j\mid j\leq b_i,\ j\ne 
w_1,\dots,w_{i-1}\}$, for $i=2,\dots,n$; 
cf.\ Proposition~\ref{prop:132-avoiding}.

\begin{theorem}
Let $w\in S_n^{312}$ be a $312$-avoiding permutation.  
Let $b=b(w)$ be the corresponding element of $\tilde\CC_n$.
Let $\lambda=(\lambda_1,\dots,\lambda_n)$ be a partition.
Then the Demazure character $ch_{\lambda,w}$ 
equals the flagged Schur polynomial:
$$
ch_{\lambda,w}(z_1,\dots,z_n) = s_\lambda^b(z_1,\dots,z_n).
$$
\label{th:flaggedSchur=Demazure}
\end{theorem}

This theorem follows from a general result by Reiner and 
Shimozono~\cite{RS}, who expressed any flagged skew Schur polynomial
as a combination of Demazure characters (key polynomials).
Theorem~\ref{th:flaggedSchur=Demazure}
implies that every Schubert polynomial $\S_w$, 
for a vexillary permutation $w\in S_n^{2143}$,
is equal to some Demazure character $ch_{\lambda,u}$, for a certain
$312$-avoiding permutation $u\in S_m^{312}$, $m<n$, associated with $w$.
Let us give a simple proof of Theorem~\ref{th:flaggedSchur=Demazure}.

Let $b=(b_1,\dots,b_n)\in\tilde\CC_n$.  Let us say that $k\in\{1,\dots,n-1\}$
is an {\it isolated entry\/} in $b$ if $k$ appears in the sequence $b$
exactly once.  Let us write $b \,{\buildrel k\over\longrightarrow}\, b'$
if $k$ is an isolated entry in $b$ and $b'\in\tilde\CC_n$ is obtained from $b$ 
by adding 1 to this entry.
In other words, we have $b_{i-1}<b_i=k< b_{i+1}$, for some 
$i\in\{1,\dots,n-1\}$ (assuming that $b_0=0$), 
and $b'=(b_1,\dots,b_{i-1},b_i+1,b_{i+1},\dots,b_n)$.

\begin{lemma}  If $b\,{\buildrel k\over\longrightarrow}\, b'$, then
$T_k\cdot s_\lambda^{b}(z_1,\dots,z_n) = s_\lambda^{b'}(z_1,\dots,z_n)$.
\label{lem:Tkslambda}
\end{lemma}

\begin{proof}  The claim follows from the formula $s_\lambda^b = 
\det\left(h_{\lambda_i-i+j}(z_1,\dots,z_{b_i})\right)_{i,j=1}^n$,
the fact that the operator $T_k$ commutes with multiplication by 
$h_m(x_1,\dots,x_l)$ for $k\ne l$; and $T_k\cdot h_m(x_1,\dots,x_k)
=h_m(x_1,\dots,x_{k+1})$.
\end{proof}

Let us also write $w\,{\buildrel k\over\longrightarrow}\, w'$, for
$w,w'\in S_n$, if $w' = s_k\,w$ and $\ell(w')=\ell(w)+1$.

\begin{lemma}  For $w,w'\in S_n^{312}$, 
if $b(w)\,{\buildrel k\over\longrightarrow}\, b(w')$ then
$w\,{\buildrel k\over\longrightarrow}\, w'$.
\label{lem:bb'-ww'}
\end{lemma}

\begin{proof}
Assume $b(w)=b$, $b(w')=b'$, and
$b\,{\buildrel k\over\longrightarrow}\, b'$.
Let $b_i=k$ be the isolated entry in $b$ that we increase.  
The construction
of the map $b\mapsto w$ implies that $w_i = k$ and $w_j=k+1$ for 
some $j>i$.  It also implies that $b'\mapsto s_k w$.
The permutation $s_k w$ is obtained from $w$ by switching 
$w_i$ and $w_j$, and its length is $\ell(w)+1$.
Thus $w\,{\buildrel k\over\longrightarrow}\, w'$.
\end{proof}

\begin{exercise}  
Check that 
$b(w)\,{\buildrel k\over\longrightarrow}\, b(w')$
if and only if
$w\,{\buildrel k\over\longrightarrow}\, w'$. 
\end{exercise}


%
%
%
%
%
%
\begin{proof}[Proof of Theorem~\ref{th:flaggedSchur=Demazure}]
Let $b=b(w)\in\tilde\CC_n$.  
We claim that there is a directed path 
$b^{(0)}\,{\buildrel k_1\over\longrightarrow}\, b^{(1)}
{\buildrel k_2\over\longrightarrow}\,
\cdots \,{\buildrel k_l\over\longrightarrow}\, b^{(l)}$
from $b^{(0)}=(1,\dots,n)$ to $b^{(l)}=b$.
In other words, we can obtain the sequence $b$ from the sequence $(1,\dots,n)$
by repeatedly adding 1's to some isolated entries.   
One possible choice of such a path is given by the following rule.
We have $b_n=n$.  Let us first increase the $(n-1)$-st entry until we obtain 
$b_{n-1}$; then increase the $(n-2)$-nd entry until we obtain $b_{n-2}$, etc.

For example,
for the sequence $b=(3,3,3,5,5)$ that corresponds to $w=32154$, 
we obtain the path
$$
(1,2,3,4,5)
\,{\buildrel 4 \over\longrightarrow}\,
(1,2,3,5,5)
\,{\buildrel 2 \over\longrightarrow}\,
(1,3,3,5,5)
\,{\buildrel 1 \over\longrightarrow}\,
(2,3,3,5,5)
\,{\buildrel 2 \over\longrightarrow}\,
(3,3,3,5,5).
$$
This path gives the reduced decomposition $s_2 s_1 s_2 s_4$
for $w=32154$.

If $w=\id$ then $b(w)=(1,\dots,n)$ and 
$\ch_{\id,\lambda}=s_\lambda^{(1,\dots,n)}=
z_1^{\lambda_1}\dots z_n^{\lambda_n}$.  
In general, according to Lemmas~\ref{lem:Tkslambda} and~\ref{lem:bb'-ww'}, 
we have $w= s_{k_l}\cdots s_{k_1}$, and thus, 
$s_\lambda^b = T_w(s_\lambda^{(1,\dots,n)}) = T_w (x^\lambda)=\ch_{w,\lambda}$.
%
\end{proof}

\begin{remark} 
Lemma~\ref{lem:bb'-ww'}, together with the exercise,
gives a bijective correspondence between
paths $(1,\dots,n) \,{\buildrel k_1 \over\longrightarrow}\,\cdots
\,{\buildrel k_l \over\longrightarrow}\,b(w)$ and the special class
of reduced decompositions $w = s_{k_l}\cdots s_{k_1}$ such that
all truncated decompositions $s_{k_i}\cdots s_{k_1}$ give 312-avoiding
permutations, for $i=1,\dots,l$.
%
\end{remark}

\begin{corollary}
Let us use the notation of Theorem~\ref{th:flaggedSchur=Demazure}.
The dimension of the Demazure module is given by the following 
matrix of binomial coefficients:
$$
\dim V_{\lambda,w} = 
\det\left(\binom{\lambda_i+b_i-i}{b_i-j}\right)_
{i,j=1}^n .
$$
\label{cor:dimVlambdaw}
\end{corollary}

\begin{proof}
We have $\dim V_{\lambda,w} = \ch_{\lambda,w}(1,\dots,1)$.
The claim follows from the determinant 
expression~(\ref{eq:flagged-schur=det}) for 
the flagged Schur polynomial 
$\ch_{\lambda,w}=s_\lambda^{(1,\dots,n),b}$ 
and the fact that $h_m^{[k,l]}(1,\dots,1) = \binom{l-k+m}{l-k}$.
\end{proof}

Corollary~\ref{cor:dimVlambdaw} presents $\dim V_{\lambda,w}$
as a polynomial of degree $\sum (b_i-i)=\ell(w)$.
According to Proposition~\ref{prop:demazure-limit},
the leading homogeneous component of this polynomial equals
$\D_w(\lambda)$.
Thus Corollary~\ref{cor:dimVlambdaw} produces
the same determinant expression
$\D_w(\lambda)=\det\left(\lambda_i^{(b_i-j)}\right)$
for a 312-avoiding permutation $w$ as Theorem~\ref{th:3412-avoiding}.

Let us give another expression for the Demazure characters 
$\ch_{\lambda,w}$ that generalizes the Weyl character formula.  
It is not hard to prove it by induction similar
to the above argument.

\begin{proposition}
Let $w\in S_n^{312}$ be a 312-avoiding permutation and let $b(w)=
(b_1,\dots,b_n)$.  
Let $W_b=\{u\in S_n\mid u_i\leq b_i,\textrm{ for any } i=1,\dots,n\}$,
and let $\Phi_{u,b}^+=
\{\epsilon_i-\epsilon_j\mid 1\leq i<j\leq b_{u^{-1}(i)}\}\subseteq\Phi^+$.
Then
$$
\ch_{\lambda,w}(z_1,\dots,z_n) = \sum_{u\in W_b} (-1)^{\ell(u)} \,
z^{u(\lambda+\rho)-\rho}
\, \prod_{\alpha\in\Phi_{u,b}^+}(1-z^{-\alpha})^{-1}.
$$
\end{proposition}

The set $W_b$ is in one-to-one correspondence 
with rook placements in the Young diagram of shape $(b_n,b_{n-1},\dots,b_1)$.
We have $|W_b| = b_1\,(b_2-1)\,(b_3-2)\,\cdots (b_n-n+1)$.
For any $u\in W_b$, we have $|\Phi_{u,b}^+|=\ell(w)$.

\section{Generalized Gelfand-Tsetlin polytope}
\label{sec:generalized-GT}

In this section we show how flagged Schur functions and Demazure
characters are related
to generalized Gelfand-Tsetlin polytopes studied by Kogan~\cite{Kog}.
\smallskip

A {\it Gelfand-Tsetlin pattern\/} $P$ of size $n$ is a triangular array
of real numbers $P=(p_{ij})_{n\geq i\geq j\geq 1}$ that satisfy 
the inequalities $p_{i-1\,j-1}\geq p_{ij}\geq p_{i-1\,j}$.
These patterns are usually arranged on the plane as follows:
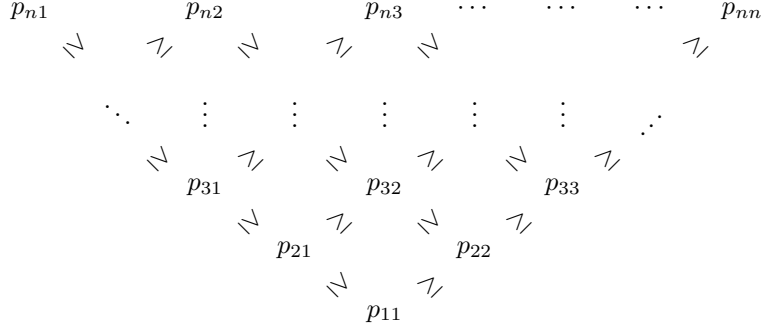
\def\Dgeq{\pspicture(0,0)\rput{-45}(0,0){$\geq$}\endpspicture}
\def\Ugeq{\pspicture(0,0)\rput{45}(0,0){$\geq$}\endpspicture}
$$
\begin{array}{ccccccccccccccccc}
p_{n1}&&&& p_{n2} &&&& p_{n3} &&\cdots&& \cdots &&\cdots&&  p_{nn}\\
&\Dgeq&&\Ugeq&&\Dgeq&&\Ugeq&&\Dgeq&&&&&&\Ugeq&\\
&&&&&&&&&&&&&&&&\\
&&\ddots&&\vdots&&\vdots&&\vdots &&\vdots&&\vdots&&\rput{45}(.1,.2){\cdots}&&\\
&&&\Dgeq&&\Ugeq&&\Dgeq&&\Ugeq&&\Dgeq&&\Ugeq&&&\\
&&&&p_{31}&&&&p_{32}&&&&p_{33}&&&&\\
&&&&&\Dgeq&&\Ugeq&&\Dgeq&&\Ugeq&&&&&\\
&&&&&&p_{21}&&&& p_{22}&&&&&&\\
&&&&&&&\Dgeq&&\Ugeq &&&&&&&\\
&&&&&&&&p_{11}&&&&&&&&\\
\end{array}
$$
The {\it shape\/} $\lambda=(\lambda_1,\dots,\lambda_n)$
of a Gelfand-Tsetlin pattern $P$ is given by $\lambda_i=p_{ni}$, 
for $i=1,\dots,n$, i.e., the shape is the top row of a pattern.
The {\it weight\/} $\beta=(\beta_1,\dots,\beta_n)$ of a Gelfand-Tsetlin 
pattern $P$ is given by $\beta_1=p_{11}$ and
$\beta_i = p_{i1}+\dots p_{ii} - p_{i-1\, 1}-\dots - p_{i-1\,i-1}$,
for $i=2,\dots,n$, i.e., the $i$-th row sum $p_{i1}+\cdots+p_{ii}$ 
equals $\beta_1+\cdots+\beta_i$.

The {\it Gelfand-Tsetlin polytope\/} $\P_\lambda\in\R^{\binom n2}$
is the set of all Gelfand-Tsetlin patterns of shape $\lambda$.  
This is a a convex polytope.  
A Gelfand-Tsetlin pattern $P=(p_{ij})$ is called {\it integer\/} 
if all $p_{ij}$ are integers.  The integer Gelfand-Tsetlin patterns 
are the lattice points of the polytope $\P_\lambda$. 

The integer Gelfand-Tsetlin patterns $P=(p_{ij})$ of shape $\lambda$ 
and weight $\beta$ are in one-to-one correspondence with semistandard 
Young tableaux $T=(t_{ij})$ of shape $\lambda$ and weight $\beta$.
 This correspondence is given by 
setting $p_{ij} = \#\{k\mid t_{kj}\leq i\}$,
i.e., $p_{ij}$ is the number of entries less than or equal to $i$
in the $j$-th row of $T$.
The proof of the following claim is immediate from the definitions.

\begin{lemma}  A semistandard Young tableau $T$ is a flagged
tableau with flags $(1,\dots,1)$ and $(b_1,\dots,b_n)$ if and only 
if the corresponding Gelfand-Tsetlin pattern $P=(p_{ij})$ satisfies
the conditions $p_{ni}=p_{n-1\,i}=\cdots= p_{b_i\,i}$, 
for $i=1,\dots,n$.
\end{lemma}

Let $w\in S_n^{312}$ be a 312-avoiding permutation, let 
$b=(b_1,\dots,b_n)=b(w)\in\tilde\CC_n$,
and let $\lambda=(\lambda_1,\dots,\lambda_n)$ be a partition.
Let us define the {\it generalized Gelfand-Tsetlin polytope\/}
$\P_{\lambda,w}$ as the set of all Gelfand-Tsetlin patterns $P=(p_{ij})$ 
of size $n$ such that 
$\lambda_i=p_{ni}=p_{n-1\,i}=\cdots=p_{b_i \,i}$, 
for $i=1,\dots,n$.  Note that $b_1+\cdots+b_n-\binom {n+1}2=\ell(w)$
is the number of unspecified entries in a pattern.
Thus $\P_{\lambda,w}$ is a convex polytope naturally 
embedded into $\R^{\ell(w)}$.  These polytopes were studied 
by Kogan~\cite{Kog}.

According to Theorem~\ref{th:flaggedSchur=Demazure}, the Demazure
character $\ch_{\lambda,w}$, for a 312-avoiding permutation $w$, is given
by counting lattice points of the generalized Gelfand-Tsetlin polytope
$\P_{\lambda,w}$.

\begin{corollary} For $w\in S_n^{312}$ and a partition 
$\lambda=(\lambda_1,\dots,\lambda_n)$, we have
$$
\ch_{\lambda,w}(z_1,\dots,z_n) = 
s_\lambda^b(z_1,\dots,z_n) =
\sum_{P\in \P_{\lambda,w}\cap \Z^{\ell(w)}} z^P,
$$
where the sum is over lattice points in the polytope $\P_{\lambda,w}$,
$z^P = z_1^{\beta_1}\cdots z_n^{\beta_n}$, and
$\beta=(\beta_1,\dots,\beta_n)$ is the weight of $P$.
In particular, the dimension of the Demazure module 
$V_{\lambda,w}$ is equal to the number of lattice points in 
the polytope $\P_{\lambda,w}$:
$$
\dim V_{\lambda,w} = \# (\P_{\lambda,w} \cap \Z^{\ell(w)}).
$$
Finally, the $\lambda$-degree of the Schubert variety $X_w$ 
divided by $\ell(w)!$ equals the  volume of the generalized Gelfand-Tsetlin
polytope $\P_{\lambda,w}$:
$$
\frac 1 {\ell(w)!} \,
\deg_\lambda(X_w) = \D_w(\lambda)=\Vol(\P_{\lambda,w}),
$$
where $\Vol$ denotes the usual volume form on $\R^{\ell(w)}$ such that 
the volume of the unit $\ell(w)$-hypercube equals $1$.
\label{cor:ch312}
\end{corollary}

The following claim is also straightforward from the definition of
the polytopes $\P_{\lambda,w}$.

\begin{proposition}  The polytope $\P_{\lambda,w}$ is the Minkowski sum of 
the polytopes $\P_{\omega_i,w}$ for the fundamental weights:
$$
\P_{\lambda,w} = a_1 \P_{\omega_1,w}+\cdots + a_{n-1} \P_{\omega_{n-1},w},
$$
where $\lambda=a_1\omega_1+\cdots+a_{n-1}\omega_{n-1}$.
\label{prop:Minkowsky}
\end{proposition}
The last claim implies that $\dim V_{\lambda,w}$ is the mixed 
lattice point enumerator of the polytopes $\P_{\omega_i,w}$, $i=1,\dots,n-1$.

\begin{remark}  
Toric degenerations of Schubert varieties $X_w$ for Kempf elements
(312-avoiding permutations in our terminology), were constructed
by Gonciulea and Lakshmibai~\cite{GL}, and were studied by
Kogan~\cite{Kog} and Kogan-Miller~\cite{KM}.
According to~\cite{Kog,KM}, these toric degenerations are associated
with generalized Gelfand-Tsetlin polytopes $\P_{\lambda,w}$. 
It is a standard fact that the degree of a toric variety is equal to the 
normalized volume of the corresponding polytope.
\end{remark}

\begin{remark}
We can extend the definition of generalized Gelfand-Tsetlin polytopes 
$\P_{w,\lambda}$ to a larger class of permutations, as follows.  
For a 231-avoiding permutation $w$, define
$\P_{w,\lambda} = \P_{w_\circ w w_\circ,\,(-\lambda_n,\dots,-\lambda_1)}$,
cf.\ Lemma~\ref{lem:involution-omega}(1).
Let $w = w^{1}\times \cdots \times w^{k}\in S_{n_1}\times \cdots \times S_{n_k}\subset S_n$
be a permutation such that all blocks $w^{i}\in S_{n_i}$ are
either 312-avoiding or 231-avoiding,  and let $\lambda$ be the
concatenation of partitions $\lambda^{1},\dots,\lambda^{k}$
of lengths $n_1,\dots,n_k$.
We have $ch(V_{\lambda,w}) = \prod ch (V_{\lambda^{i},\,w^{i}})$
and $\D_w(\lambda) = \prod \D_{w^{i}}(\lambda^{i})$.
Let us define $\P_{w,\lambda} = \P_{w^{1},\,\lambda^{1}}\times \cdots
\times\P_{w^{k},\,\lambda^{k}}$.
Then Corollary~\ref{cor:ch312} and Proposition~\ref{prop:Minkowsky} 
remain valid for this more general class of permutations
with 312- or 231-avoiding blocks.
These claims extend results of Dehy and Yu~\cite{DY}.
\end{remark}

%
%

\section{A conjectured value of $\D_w$}
\label{sec:conj}

In this section we give a conjectured value of $\D_w$ for a special
class of permutations $w$.
\smallskip

Let $w$ be a permutation whose code has the form
 $$ \code(w)=(n, *, n-1, *, n-2, \cdots, *, 2, *, 1,0,0,\dots), $$
where each $*$ is either 0 or empty. We
call such a permutation \emph{special}. For instance, $w=761829543$ is
special, with $\code(w)=(6,5,0,4,0,3,2,1,0,\dots)$. Note also that
$w_\circ$ is special.  Suppose that $w$ is special with
$\code(w)=(c_1,c_2,\dots)$. Let $c_1=n$, and let $k$ be the number of
0's in $\code(w)$ that are preceded by a nonzero number, i.e, $c_i$=0,
$c_{i-1}>0$. Let $a_1<\cdots<a_k=n+k$ be the positions of these 0's,
so $c_{a_1} = \cdots = c_{a_k}=0$. Define
 \begin{eqnarray*} a_\delta(y_1,\dots,y_n) & = &
    \prod_{1\leq i<j\leq n}(y_i-y_j)\\ & = &
   \sum_{w\in S_n} (-1)^{\ell(w)}\, y_1^{w(1)-1}\cdots y_n^{w(n)-1}.
  \end{eqnarray*}
An $n$-element subset $J=\{j_1,\dots,j_n\}$ of
$\{1,2,\dots,n+k\}$ is said to be \emph{valid} (with respect to $w$)
if
 $$ \#(J\cap \{a_{i-1}+1,a_{i-1}+2,\dots,a_i\}) = a_i-a_{i-1}-1 $$
for $1\leq i\leq k$ (where we set $a_0=0$). For instance if
$\code(w)=(3,0,2,1,0)$, then the valid sets are 134, 135, 145,
234, 235, 245. Clearly the number of valid sets in general is equal to
$(a_1-1)(a_2-a_1-1)\cdots(a_k-a_{k-1}-1)$. If $J$ is a valid set, then
define the \emph{sign} $\varepsilon_J$ of $J$ by $\varepsilon_J =
(-1)^{d_J}$, where
 $$ d_J = \binom{n+k+1}{2}-1-(a_1+1)-\cdots-(a_{k-1}+1)-
        \sum_{i\in J} i. $$
Note that the quantity $\binom{n+k+1}{2}-1-(a_1+1)-\cdots-
(a_{k-1}+1)$ appearing above is just $\sum_{i\in L}i$ for the valid
subset $L$ with largest element sum, viz.,
 $$ L=\{ 1,2,\dots,n+k\}-\{1,a_1+1,a_2+1,\dots,a_{k-1}+1\}. $$
In particular, $d_L=0$ and $\varepsilon_L=1$.

\begin{conjecture} \label{conj:dw}
Let $w$ be special as above. Then
 $$ \D_w = C_{nk}
    \sum_{J=\{j_1,\dots,j_k\}} \varepsilon_J \,
   a_\delta(y_{n+k-j_1+1},y_{n+k-j_2+1},\dots,y_{n+k-j_k+1}), $$
where
 $$ C_{nk} = \frac{(n+1)!\,(n+2)!\cdots(n+k-1)!}{\binom{n+1}{2}!} $$
and $J$ ranges over all valid subsets of $\{1,2,\dots,n+k\}$.
\end{conjecture}

As an example of Conjecture~\ref{conj:dw}, let $w=41532$, so
$\code(w)= (3,0,2,1,0)$. Write $y_1=a$, $y_2=b$, etc. Then
 $$\D_w = \frac{1}{30}(a_\delta(a,b,d)-a_\delta(a,b,e) -
  a_\delta(a,c,d)+a_\delta(a,c,e)+a_\delta(b,c,d)-
  a_\delta(b,c,e)). $$
We have verified Conjecture~\ref{conj:dw} for $n\leq 5$.

\section{Schubert-Kostka matrix and its inverse}
\label{sec:Kostka}

In this section we discuss the following three equivalent problems:
\begin{enumerate}
\item Express the polynomials $\D_w$ as linear combinations
of monomials. 
\item Express monomials as linear combinations of Schubert
polynomials $\S_w$. 
\item Express  Schubert polynomials as linear combination
of standard elementary monomials $e_{a_1}(x_1) e_{a_2}(x_1,x_2)
e_{a_3}(x_1,x_2,x_3)\cdots$.
\end{enumerate}

Let $\N^\infty$ be the set of ``infinite compositions'' 
$a=(a_1,a_2,\dots)$ such that all $a_i\in\N=\Z_{\geq 0}$ and $a_i=0$, for
almost all $i$'s.  
For $a\in \N^\infty$, let $x^a=x_1^{a_1} x_2^{a_2}\cdots$ and 
$y^{(a)} = \frac{y_1^{a_1}}{a_1!}
\frac{y_2^{a_2}}{a_2!}\cdots$.
The polynomial ring $\Q[x_1,x_2,\dots]$ in infinitely many variables
has the linear bases $\{x^a\}_{a\in \N^\infty}$ and 
$\{\S_w\}_{w\in S_\infty}$;
also the polynomial ring $\Q[y_1,y_2,\dots]$ has the linear basis
$\{y^{(a)}\}_{a\in \N^\infty}$ and $\{\D_w\}_{w\in S_\infty}$,
where $\S_w=\S_w(x_1,x_2,\dots)$ and $\D_w=\D_w(y_1,y_2,\dots)$.

Let us define the {\it Schubert-Kostka matrix\/} $K=(K_{w,a})$, 
$w\in S_\infty$ and $a\in \N^\infty$, by
$$
\S_{w} = \sum_{a\in \N^\infty} K_{w,a}\, x^a.
$$
The numbers $K_{w,a}$ are nonnegative integers.
They can be combinatorially interpreted in terms of 
{\it RC-graphs;} see~\cite{FK} and~\cite{BJS}.  
For grassmannian permutations $w$, the numbers $K_{w,a}$ are 
equal to the usual Kostka numbers, which are the coefficients of
monomials in Schur polynomials.


The matrix $K$ is invertible, because every monomial $x^a$ can be
expressed as a finite linear combination of Schubert polynomials.  Let 
$K^{-1}=(K^{-1}_{a,w})$ be the inverse of the Schubert-Kostka matrix.  We have
$$
x^a = \sum_{w\in S_\infty} K_{a,w}^{-1}\, \S_w.
$$

The basis $\{x^a\}_{a\in \N^\infty}$ is D-dual to 
$\{y^{(a)}\}_{a\in \N^\infty}$,
and the basis $\{\S_w\}_{w\in S_\infty}$ is D-dual to
$\{\D_w\}_{w\in S_\infty}$; see Corollary~\ref{cor:D-dual-sigma}.
Thus the previous two formulas are equivalent to the following statement.

\begin{proposition}  
We have
$$
y^{(a)} = \sum_{w\in S_\infty} K_{w,a} \,\D_w
\quad\textrm{and, equivalently,}\quad
\D_w = \sum_{a\in \N^\infty} K_{a,w}^{-1}\,y^{(a)}.
$$
\label{prop:DKl}
\end{proposition}

This claim shows that an explicit expression
for the polynomials $\D_w$ in terms of monomials
is equivalent to a formula for 
entries of the inverse Schubert-Kostka matrix $K^{-1}$.
We remark that a combinatorial interpretation of the inverse 
of the usual Kostka matrix was given by Egecioglu and 
Remmel~\cite{ER}.  It would be interesting to give a
subtraction-free combinatorial interpretation for entries of 
the inverse of the Schubert-Kostka matrix.    
Notice that the matrix $K^{-1}$ has both positive and negative entries.
Although we do not know such a formula in general,
it is not hard to give an alternating formula for the entries of $K^{-1}$,
as follows.

Let us fix a positive integer $n$.
Let $w_\circ$ be the longest permutation in $S_n$,  let 
$\N^n$ be the set of compositions 
$a=(a_1,\dots,a_n)$, $a_i\in\N$, naturally embedded into $\N^\infty$,
and let 
$\rho=(n-1,n-2,\dots,0)\in\N^n$.

\begin{lemma}  If $w\in S_n$, then $K_{a,w}^{-1} = 0$, unless
$a\in\N^n$.
\label{lem:K=0}
\end{lemma}

\begin{proof}  Follows from Proposition~\ref{prop:DKl}
and the fact that $\D_w$ involves
only $y_1,\dots,y_n$, for $w\in S_n$. 
\end{proof}

Assume by convention that $K_{w,a}=0$ 
if some entries $a_i$ are negative.

\begin{proposition}
Assume that $w\in S_n$.   Then, for any $a\in \N^n$, we have
$$
K_{a,w}^{-1} = \sum_{u\in S_n} (-1)^{\ell(u)}\, K_{w_\circ w, u(\rho)-a}.
$$
\label{prop:K-inv}
\end{proposition}

\begin{proof}
Follows from Corollary~\ref{cor:D-S}(2), for $u=\id$,
and Proposition~\ref{prop:DKl}.
\end{proof}

For a $312$-avoiding permutation $w$, Proposition~\ref{prop:K-inv} 
implies a more explicit expression for $K_{a,w}^{-1}$.
Indeed, in this case, $\S_{w_\circ w}  = x^{c}$,
where $c=\code(w_\circ w)$.  In other words, $K_{w_\circ w,u(\rho)-a}$
equals $1$, if $u(\rho)-a=c$, and $0$, otherwise.
We obtain the following result.

%

\begin{corollary}
For a $312$-avoiding permutation $w\in S^{312}_n$ with
$c=\code(w_\circ w)$,
and an arbitrary $a=(a_1,\dots,a_n)\in \N^n$, 
we have
$$
K_{a,w}^{-1} =
\left\{
\begin{array}{cl}
(-1)^{\ell(u)} &
\textrm{if } 
a+c = u(\rho),
\textrm{ for some permutation } u\in S_n,\\[.1in]
0 &\textrm{otherwise.}
\end{array}
\right.
$$
\label{cor:K-inv-312}
\end{corollary}

Note that this expression for $K_{a,w}^{-1}$ is stable
under the embedding $S_n\hookrightarrow S_{n+1}$.
More generally, we can give an expression for $K_{a,w}^{-1}$,
for any $3412$-avoiding permutation $w$, as a sum over flagged 
semistandard tableaux; cf.~Theorem~\ref{th:3412-avoiding}.
Also Conjecture~\ref{conj:dw} implies a conjecture for values 
$K_{a,w}^{-1}$, for special permutations $w$, as defined
in Section~\ref{sec:conj}.

Recall that the involution $y_i\mapsto - y_{n+1-i}$ sends 
$\D_w$ to $\D_{w_\circ w w_\circ}$ 
(see Lemma~\ref{lem:involution-omega}).
If $w\in S_n$, then the second identity in
Proposition~\ref{prop:DKl} involves only terms with $a\in\N^n$.
Applying the above involution to this identity,
we deduce that the inverse Schubert-Kostka matrix has the following symmetry.

\begin{lemma} For any $w\in S_n$ and $a=(a_1,\dots,a_n)\in\N^n$, we have
$$
K_{a,w}^{-1} = (-1)^{|a|} K_{\bar a, w_\circ w \,w_\circ}^{-1},
$$
where $|a| = a_1+\cdots + a_n$ and $\bar a = (a_n,\dots,a_1)$.
\end{lemma}

\begin{remark}
The matrix $K$ does not have this kind of symmetry.  For example, 
$\S_{s_1}=x_1$ and $\S_{s_{n-1}} = x_1 + \cdots + x_{n-1} \ne - x_1$.  
Thus $K_{s_1,(1\,0^{n-1})} = 1$ and $K_{w_\circ s_1 w_\circ,(0^{n-1} 1)} = 0$.  
An argument similar to the above does not work for matrix $K$, because
the first identity in Proposition~\ref{prop:DKl} may involve terms 
with $w\in S_\infty\setminus S_n$ even if $a\in \N^n$.
\end{remark}

Applying this symmetry to Corollary~\ref{cor:K-inv-312}, we obtain
an explicit expression for $K^{-1}_{a,w}$, for 231-avoiding permutations $w$,
as well.

\begin{corollary}
For a $231$-avoiding permutation $w\in S^{231}_n$
with $\code(w\,w_\circ)=(c_1,\dots,c_n)$
and an arbitrary $a=(a_1,\dots,a_n)\in \N^n$, 
we have
$$
K_{a,w}^{-1} =
\left\{
\begin{array}{cl}
(-1)^{\ell(u)+|a|} &
\textrm{if } 
(c_1+a_n,\dots,c_n+a_1)=u(\rho),
\textrm{ for some } u\in S_n,
\\[.1in]
0 &\textrm{otherwise.}
\end{array}
\right.
$$
\label{cor:K-inv-231}
\end{corollary}

Say that a permutation $w$ is {\it strictly dominant\/}
if its code $\code(w) = (c_1,\dots,c_n)$
is a strict partition, i.e.,
$c_1>c_2>\cdots> c_k = c_{k+1} = \cdots  = c_n = 0$,
for some $k=1,\dots,n$.

\begin{exercise}  {\rm (A)} \ Show that the following conditions are 
equivalent:
\begin{enumerate}
\item $w$ is strictly dominant;
\item $w w_\circ$ is strictly dominant;
\item $w$ is of the form $w_1>w_2>\cdots > w_k<w_{k+1}<\cdots < w_n$;
\item $w$ is both 132-avoiding and 231-avoiding.
\end{enumerate}

\noindent
{\rm (B)} \ There are exactly $2^{n-1}$ strictly dominant permutations in $S_n$.

\noindent
{\rm (C)} \
If $w$ is strictly dominant with 
$\code(w) = (c_1>\cdots> c_{k-1} > 0 = \cdots  = 0)$, then
$\code(w w_\circ) = (c_1'> \cdots > c_{n-k}'>0= 
\cdots = 0)$, where the set $\{c_1',\dots,c_{n-k}'\}$ 
is the complement to the set $\{c_1,\dots,c_{k-1}\}$
in $\{1,\dots,n-1\}$.
\end{exercise}

Let us specialize Corollary~\ref{cor:K-inv-231} to strictly 
dominant permutations.

\begin{corollary}  Let $w$ be a strictly dominant permutation
with $\code(w) = (c_1>\cdots > c_{k-1} > c_k =\cdots =0)$.
Assume that $a = (a_1,\dots,a_k,0,\dots,0)$.
Then
$$
K_{a,w}^{-1} = 
\left\{
\begin{array}{cl}
(-1)^{\ell(\sigma)} &
\textrm{if } (a_1,\dots,a_k) = (c_{\sigma_1},\dots,c_{\sigma_k}), 
\textrm{ for some } \sigma\in S_k,\\[.1in]
0 &\textrm{otherwise.}
\end{array}
\right.
$$
Equivalently, we have
$\D_w(y_1,\dots,y_k,0,\dots,0) = \sum_{\sigma\in S_k}
(-1)^{\ell(\sigma)}\, y_{\sigma_1}^{(c_1)}\cdots y_{\sigma_k}^{(c_k)}$.
\label{cor:K-strictly-dominant}
\end{corollary}

\begin{proof}  We have $\code(w w_\circ) = 
(c_1'>\cdots > c_{n-k}'>0=\cdots =0)$, where
$\{c_1',\dots,c_{n-k}'\}$ is the set complement 
$\{0,\dots,n-1\}\setminus \{c_1,\dots,c_{k}\}$.
According to Corollary~\ref{cor:K-inv-231}, $K_{a,w}^{-1}=0$, unless
$c_1',\dots,c_{n-k}',a_k,\dots,a_1$ is a permutation of $0,\dots,n-1$;
or, equivalently, $a_1,\dots,a_k$ is a permutation of $c_1,\dots,c_k$.
We leave it as an exercise for the reader to check that the signs agree.
\end{proof}

According to Lemma~\ref{lem:K=0},
for the strictly dominant permutation 
$w=(k,k-1,\dots,1,k+1,k+2,\dots,n)\in S_k\subset S_n$,
the assertion of Corollary~\ref{cor:K-strictly-dominant}
is true for an arbitrary $a$, without the assumption that
$a=(a_1,\dots,a_k,0,\dots,0)$.   
However, if we skip this assumption, for other permutations, 
we will have more cases.
For example, for $w=(k+1,k-1,\dots,1,k,k+2,k+3,\dots,n)$
with $\code(w) = (k,k-2,\dots,1,0,\dots,0)$, 
Corollary~\ref{cor:K-inv-231}  implies that
$$
K_{a,w}^{-1} = 
\left\{
\begin{array}{cl}
(-1)^{\ell(\sigma)} &
\textrm{if } a = (c_{\sigma_1},\dots,c_{\sigma_k},0,\dots,0),
\textrm{ for some } \sigma\in S_k,\\[.1in]
(-1)^{\ell(\tau)+1} &
\textrm{if } a = (k-\tau_1,\dots,k-\tau_k,1,\dots,0), 
\textrm{ for some } \tau\in S_k,\\[.1in]
0 &\textrm{otherwise.}
\end{array}
\right.
$$

The polynomial ring $\Q[x_1,x_2,\dots]$ has the following
basis of {\it standard elementary monomials:}  
$e_{a} := e_{a_2}(x_1)\, e_{a_3}(x_1,x_2) \, e_{a_4}(x_1,x_2,x_3)\cdots$, 
where $a=(a_1,a_2,\dots)\in\N^\infty$ such that $0\leq a_i\leq i-1$,  
for $i=1,2,\dots$.
This basis was originally introduced by Lascoux and 
Sch\"utzenberger~\cite{LS}; see also~\cite[Proposition~3.3]{FGP}.

\begin{remark} Expressions for Schubert polynomials 
in the basis of standard elementary monomials play
an important role in calculation of Gromov-Witten invariants for
the small quantum cohomology ring of the flag manifold; see~\cite{FGP}.
\end{remark}

The Cauchy formula (Lascoux~\cite{La1}, see also, e.g.,~\cite{Man})
$$
\sum_{w\in S_n} \S_w(x)\cdot \S_{w w_\circ} (y)=
\prod_{i+j\leq n}(x_i+y_j) = \prod_{k=1}^{n-1}\sum_{i=0}^k y_{n-k}^{k-i}
e_i(x_1,\dots,x_k)
$$
implies that
$$
e_{w_\circ(\rho-a)} = \sum_{w\in S_n} K_{w,a}\, \S_{w w_\circ},
$$
for $a\in\N^n$.
Equivalently,
$$
\S_{w w_\circ} = \sum_a K^{-1}_{a,w}\, e_{w_\circ(\rho-a)}.
$$
This shows that the problem of inverting the Schubert-Kostka matrix
is equivalent to the problem of expressing a Schubert polynomial
in the basis of standard elementary monomials.

Let us assume, by convention, that $e_a=0$, unless $0\leq a_i\leq i-1$,
for $i\geq 1$.  Proposition~\ref{prop:K-inv} implies the following claim.

\begin{corollary}  For $w\in S_n$, the Schubert polynomial $\S_w$ 
can be expressed as
$$
\S_w = \sum_{u\in S_n,\,a\in\N^n} (-1)^{\ell(u)} K_{w_\circ w w_\circ,\,
w_\circ(a) + u(\rho) -\rho}\, e_a.
$$
\end{corollary}

In particular, for 213-avoiding permutations, we obtain the following 
result.

\begin{corollary}
For a $213$-avoiding permutation $w\in S_n$
and $c=\code(w_\circ w w_\circ)$, the Schubert $\S_w$
polynomial can be expressed as
$$
\S_w = \sum_{u\in S_{n-1}} (-1)^{\ell(u)}\, e_{w_\circ(c+\rho-u(\rho))}.
$$
\end{corollary}

Let us also give a (not very difficult) alternating expression for 
the generalized Littlewood-Richardson coefficients.

\begin{corollary}  Let $u,v,w\in S_n$.  Then the generalized
Littlewood-Richardson coefficient $c_{u,v,w}$ is equal to 
$$
c_{u,v,w} = 
\sum_{a,b} K_{u,a}\,K_{v,b}\,K_{a+b,w_\circ w}^{-1}
=
\sum_{z,a,b,c} (-1)^{\ell(z)} K_{u,a}\,K_{v,b}\,K_{w,c}, 
$$
where the second sum is over permutations $z\in S_n$ 
and compositions $a,b,c\in \N^n$ such that $a+b+c = z(\rho)$.
\end{corollary}

\begin{proof}
We have 
$\S_u\cdot \S_v = \sum_{a,b} K_{u,a}\,K_{v,b}\,x^{a+b}
= \sum_{a,b,w} K_{u,a}\,K_{v,b}\,K_{a+b,w_\circ w}^{-1}\,\S_{w_\circ w}$,
which implies the first claim.
Now apply Proposition~\ref{prop:K-inv}.
\end{proof}


Let us identify the polynomial rings 
$\Q[x_1,x_2,\dots]=\Q[y_1,y_2,\dots]$.
The transition matrix between the bases $\{\S_w\}$ and $\{x^a\}$
is $K$;
the transition matrix between the bases $\{x^a\}$ and $\{x^{(a)}\}$
is the diagonal matrix $D$ with products of factorials;
and the transition matrix between the bases $\{x^{(a)}\}$ and $\{\D_w\}$
is $K^T$.  Thus the transition matrix between the bases
$\{\S_w\}$ and $\{\D_u\}$ is $KDK^T$.  In other words,
we obtain the following result.

\begin{corollary}
We have $\S_u = \sum_{w\in S_\infty} L_{u,w}\,\D_w$, where
$$
L_{u,w} = \sum_{a\in \N^\infty} K_{u,a}\, \,K_{w,a}\,a_1!\, a_2!\cdots
=(\S_u,\S_u)_D.
$$
\end{corollary}

Notice that the matrix $L$ is symmetric, i.e.,
the coefficient of $\D_w$ in $\S_u$ equals the coefficient of
$\D_u$ in $\S_w$.


\section{Parking functions}
\label{sec:parking}

Let $n=r+1$.
Assume that $w=(1,2,\dots,r+1)=s_1 s_2\cdots s_r\in S_{r+1}$ 
is the long cycle.  In this section we calculate the corresponding 
polynomial $\D_r = \D_{s_1\dots s_r}$ in five different ways.
\smallskip

Let us use the coordinates $Y_i=(y,\alpha_i^\vee)$, $i=1,\dots,r$, from 
Section~\ref{sec:examp-duan}.  These coordinates
are related to the coordinates $y_1,\dots,y_{r+1}$ from
Section~\ref{sec:degree-type-A} by
$Y_i=y_i-y_{i+1}$, for $i=1,\dots,r$.
In the notation of Corollary~\ref{cor:duan}, for $w=s_1\cdots s_r$, 
we have $(i_1,\dots,i_l)=(1,\dots,r)$, and the Cartan integer 
$a_{i_p i_q}$ is $-1$, if $q=p+1$, and $0$, if $q>p+1$.  
Thus the sum in Corollary~\ref{cor:duan} involves only terms 
corresponding to arrays $(k_{pq})$ with $k_{pq}=0$, unless $q= p+1$.  
In this case, the product $\prod k_{pq}!$ cancels with
the product $\prod K_{*s}!$.  
More explicitly, Corollary~\ref{cor:duan} gives 
$$
\D_r = \sum_{c_1,\dots,c_r}  \frac{Y_1^{c_1}}{c_1!}\cdots
\frac{Y_r^{c_r}}{c_r!},
$$
where the sum is over nonnegative integer sequence $(c_1,\dots,c_r)$
such that $c_1\leq 1,c_1+c_2\leq 2,c_1+c_2+c_3\leq 3$,\dots,
$c_1+\cdots+c_{r-1}\leq r-1$,
$c_1+\cdots+c_{r}=r$.
There are exactly the Catalan number $\frac{1}{r+1}\binom {2r}{r}$ 
of such sequences.

A {\it parking function\/} of length $r$ is a sequence 
of positive integers $(b_1,\dots,b_r)$, $1\leq b_i\leq r$, such that 
$\#\{i\mid b_i \leq k\}\geq k$, for $k=1,\dots,r$.
The number of parking functions of length $r$ equals $(r+1)^{r-1}$.
Recall that the number $(r+1)^{r-1}$ also equals the number of
spanning trees in the complete graph $K_{r+1}$.
Let us define the $r$-th {\it parking polynomial} by
$$
P_r(Y_1,\dots,Y_r) = \sum_{(b_1,\dots,b_r)} Y_{b_1}\cdots Y_{b_r},
$$
where the sum is over parking functions $(b_1,\dots,b_r)$ of length $r$.
For example,
$$
P_3 = 6\,Y_1 Y_2 Y_3 + 3\,Y_1^2 Y_2 + 3\,Y_1Y_2^2  +
3\,Y_1^2Y_3 + Y_1^3.
$$
The polynomial $\frac{1}{r!} P_r(Y_1,\dots,Y_r)$
appeared in~\cite{SP} as the volume of a certain polytope;
see Corollary~\ref{cor:SP} below.
According to~\cite{SP}, for a partition $\lambda = 
(\lambda_1,\dots,\lambda_{r+1})$,
the value $P_r(\lambda_1-\lambda_2,\lambda_2-\lambda_3,\dots,
\lambda_r-\lambda_{r+1})$ equals the number of
$\lambda$-parking functions, which generalize the usual
parking functions.

We can write the above expression for $\D_{r}$
in terms of the parking polynomial.

\begin{proposition}  We have 
$\D_{r}= \frac{1}{r!}\,P_r(Y_r,\dots,Y_2,Y_1)$.
In particular, the degree of the Schubert variety 
$X_{s_1\dots s_r}$ equals the number of trees
$$
\deg(X_{s_1\cdots s_r})=P_r(1,\dots,1)=(r+1)^{r-1}.
$$
\label{prop:s1___s_r}
\end{proposition}

\begin{remark}
Proposition~\ref{prop:s1___s_r} is true for an arbitrary Weyl group $W$
and $w=s_{i_1}\cdots s_{i_r}\in W$ such that
$(\alpha_{i_p}^\vee,\alpha_{i_{p+1}}) = 1$ and 
$(\alpha_{i_p}^\vee,\alpha_{i_{q}}) = 0$, for $q>p+1$;
see Corollary~\ref{cor:duan}.
\end{remark}

\begin{remark}
Let us weight the covering relation $u\lessdot u s_{ij}$, $i<j$,
in the Bruhat order on $S_{r+1}$ by $j-i$. 
According  to Proposition~\ref{prop:s1___s_r}, the weighted sum
over saturated chains from $\id$ to $s_1\cdots s_r$ equals the number $(r+1)^{r-1}$
of trees.  Compare this with the fact that the total number of decompositions
of the cycle $s_1\cdots s_r$ into a product of $r$ transpositions also equals 
$(r+1)^{r-1}$.
\end{remark}

Let us write the polynomial $\D_r=\D_r(y_1,\dots,y_{r+1})$
in terms of the variables $y_1,\dots,y_{r+1}$.
According to Corollary~\ref{cor:D-S}(3), the polynomial $\D_r$
is recursively given by the integration $\D_r = I_r(\D_{r-1})$. In other words,
\begin{equation}
\D_r(y_1,\dots,y_{r+1}) = \int_{y_{r+1}}^{y_r}
\D_{r-1}(y_1,\dots,y_{r-1},t)\,dt.
\label{eq:D-r-int}
\end{equation}
\begin{equation}
\D_r(y_1,\dots,y_{r+1}) =  \int_{y_{r+1}}^{y_r} dt_r \int_{t_r}^{y_{r-1}}
dt_{r-1} \cdots \int_{t_3}^{y_2}dt_2 \int_{t_2}^{y_1} dt_1.
\label{eq:D-r-int-long}
\end{equation}
Equivalent integral formulas for the parking polynomials were given 
by Kung and Yan~\cite{KY}.
The right-hand side of the second formula is easily seen to be
equal to the volume of the polytope from~\cite{SP}, see below.

The long cycle $w = s_1 \cdots s_r$ is a $312$-avoiding permutation in
$S_{r+1}$.  The code of the permutation $w_\circ w$ equals $\code(w_\circ
w)=(r-1,r-2,\dots,1,0,0)$. According to Theorem~\ref{th:3412-avoiding},
the polynomial $\D_r$ is given by the determinant
of the following almost lower-triangular $(r+1)\times (r+1)$-matrix:
\begin{equation}
\D_{r}(y_1,\dots,y_{r+1}) = \det
\begin{pmatrix}
y_1      & 1   & 0 & \cdots & 0 & 0   \\[.1in]
y_2^{(2)} & y_2 & 1 & \cdots & 0 & 0 \\[.1in] 
y_3^{(3)} & y_3^{(2)} & y_3 & \cdots & 0 & 0 \\[.1in]
\vdots & \vdots & \vdots & \ddots & \vdots & \vdots \\[.1in]
y_r^{(r)} & y_r^{(r-1)} & y_r^{(r-2)} & 
\cdots & y_r & 1 \\[.1in]
y_{r+1}^{(r)} & y_{r+1}^{(r-1)} & y_{r+1}^{(r-2)} & 
\cdots & y_{r+1} & 1 \\[.1in]
\end{pmatrix}, 
\label{eq:det-expr-1}
\end{equation}
where, as before, $y_i^{(a)} = \frac{y_i^a}{a!}$.

\begin{remark} 
Determinant~(\ref{eq:det-expr-1}) is closely related to the formula 
found by Steck \cite{Steck} and Gessel~\cite{Ges} that can be written in
our notation as
\begin{equation}
\D_{r}(y_1,\dots,y_r,0) = \det\left(y_i^{(j-i+1)}\right)_{i,j=1}^r.
\label{eq:det-expr-2}
\end{equation}
Since $\D_{r}(y_1+c,\dots,y_{r+1}+c) = \D_{r}(y_1,\dots,y_{r+1})$,
expression~(\ref{eq:det-expr-2}) defines the polynomial $\D_r$.
Expression (\ref{eq:det-expr-2}) is obtained from~(\ref{eq:det-expr-1})
by setting $y_{r+1}=0$.  On the other hand, 
we can obtain expression~(\ref{eq:det-expr-1}) for 
$\D_{r-1}$ by differentiating~(\ref{eq:det-expr-2}) with respect to $y_r$.
This implies that
$$
\D_{r-1}(y_1,\dots,y_r) = \frac{\partial}{\partial y_r}
\D_{r}(y_1,\dots,y_r,0),
$$
which is equivalent to~(\ref{eq:D-r-int}).
Kung and Yan~\cite[Sect.~3]{KY} derived this expression 
in terms of Gon\u{c}arov polynomials.
\end{remark}

Expanding the determinant~(\ref{eq:det-expr-1}), we obtain the following 
result.

\begin{proposition}  We have
$$
\D_{r} = \sum 
(-1)^{r+1-k}\,y_{i_1}^{(i_1)}\,y_{i_1+i_2}^{(i_2)}\cdots
y_{i_1+\cdots + i_{k-1}}^{(i_{k-1})}
y_{i_1+\cdots + i_{k}}^{(i_k-1)},
$$
where the sum is over $2^{r+1}$ sequences $(i_1,\dots,i_k)$
such that $i_1,\dots,i_k\geq 1$ and 
$i_1+\cdots+i_k=r+1$.
(Notice that the power of the last term is decreased by $1$.)
\end{proposition}

\begin{corollary} For $a\in\N^{r+1}$, the element 
$K_{a,\,s_1\cdots s_r}^{-1}$ of the  inverse Schubert-Kostka
matrix equals $(-1)^{r+1-k}$, if the sequence $(a_1,\dots,a_r,a_{r+1}+1)$
is the  concatenation of $k$ sequences of the form $(0,\dots,0,l)$
with $l-1$ zeros, for $l \geq 1$;
otherwise $K_{a,\,s_1\cdots s_r}^{-1}=0$.
\end{corollary}

For example, we have 
$K_{(1,0,0,3,0,2,1,0,0,2),\,s_1\cdots s_9}^{-1} = (-1)^{10-5}$.


The generalized Gelfand-Tsetlin polytope $\P_{\lambda,w}$
from Section~\ref{sec:generalized-GT},
for the 312-avoiding permutation $w=s_1\cdots s_r$, is given by 
the inequalities:
$$
\P_{\lambda,s_1\cdots s_r} = \{(t_1,\dots,t_r)\in \R^r\mid 
\lambda_i\geq t_i,\textrm{ for } i =1,\dots,r; \
t_1\geq t_2\geq \cdots\geq t_r\geq \lambda_{r+1}\}.
$$
This polytope is exactly the polytope studied in~\cite{SP}.
According to Corollary~\ref{cor:ch312}, 
$\D_{r} (\lambda)$ equals the volume of the polytope 
$\P_{\lambda,s_1\cdots s_r}$.
Also, as we already mentioned, this volume equals the
right-hand side of~(\ref{eq:D-r-int-long}),
for $(y_1,\dots,y_{r+1})=(\lambda_1,\dots,\lambda_{r+1})$.
We recover the following result from~\cite{SP} about the relation
of this polytope with the parking polynomial $P_r$.

\begin{corollary} 
We have $\Vol(\P_{\lambda,s_1\cdots s_r}) = \frac{1}{r!}
P_r(Y_r,\dots,Y_1)$,
where $Y_i = \lambda_i -\lambda_{i+1}$, for $i=1,\dots,r$.
\label{cor:SP}
\end{corollary}


Let us also calculate the polynomial $\D_{r}$ using
just its definition in terms of saturated chains in the Bruhat order.

For an arbitrary Weyl group $W$ and $w=s_{i_1}\cdots s_{i_l}\in W$
with distinct $i_1,\dots,i_l$, the interval 
$[\id,w]\subset W$ in the Bruhat order consists of the elements
$u=s_{j_1}\cdots s_{j_s}$ such that $j_1,\dots,j_s$ is a subword
of $i_1,\dots,i_l$; see Section~\ref{sec:notations}.
Thus the interval $[\id,w]$ is isomorphic to the Boolean lattice of order $l$.

In particular, this is true for the long cycle
$w=s_1\cdots s_r= (1,\dots,r+1)$ in $S_{r+1}$.
The elements $u$ covered by $w$ are 
of the form $u = s_1\cdots \widehat{s_k} \cdots s_r =
w \,s_{k,r+1} = (1,2,\dots,k)(k+1,k+2,\dots,r+1)$, 
for some $k\in\{1,\dots,r\}$.
Moreover, for such $u$, the Chevalley multiplicity equals
$m(u\lessdot w) = y_k - y_{r+1}=Y_{k} + Y_{k+1} + \cdots + Y_r$.  
The interval 
$[\id,(1,\dots,k) (k+1,\dots,r+1)]$ in the Bruhat order is isomorphic to
the product of two intervals $[\id,(1,\dots,k)]\times[\id, (k+1,\dots,r+1)]$.
Thus we obtain the following recurrence relation 
for the parking polynomial $P_r$ (related to $\D_{r}$
by Proposition~\ref{prop:s1___s_r}):
$$
P_r(Y_1,\dots,Y_r) = \sum_{k=1}^r (Y_1+\cdots + Y_k)\cdot 
P_{k-1}(Y_1,\dots,Y_{k-1})\cdot P_{r-k}(Y_{k+1},\dots,Y_r).
$$
Also $P_0 = 1$ and $P_1(Y_1) = Y_1$.
This relation follows from results of Kreweras~\cite{Kre}.
It implies the following combinatorial interpretation
of the parking polynomial $P_r(Y_1,\dots,Y_r)$.

An {\it increasing binary tree\/} is a directed rooted tree with 
an increasing labeling of vertices by the integers $1,\dots,r$ such that
each vertex has at most one left successor and at most one right successor.
Let $\T_r$ be the set of such trees with $r$ vertices.  It is well known that
$|\T_r|=r!$; see~\cite{EC1}.  
Let us define the weight of a tree in $\T_r$ as follows.
For $T\in \T_r$, let $\tilde T$ be the binary tree obtained from $T$ by adding 
two leaves (left and right) to each vertex of $T$ without successors and one 
left (resp., right) leaf to each vertex of $T$ with only a right (resp., left) 
successor.
Then $\tilde T$ has $r+1$ leaves.  Let us label these
leaves by the variables $Y_1,\dots,Y_{r+1}$ from left to right.
For each vertex $v$ in $T$, define the weight $\wt(v)$ as the
sum of $Y_i$'s corresponding to all leaves of $\tilde T$ in the left branch of
$v$.  Let us define the weight of $T\in\T_r$ as  the product $\wt(T) = \prod
\wt(v)$ over all vertices $v$ of~$T$.

\psset{unit=1.4pt} 
\begin{figure}[ht]
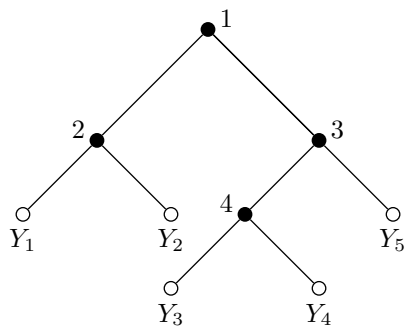

\pspicture(-70,-77)(70,7)
\cnode*(0,0){2}{A}
\rput(5,3){$1$}
\cnode*(-30,-30){2}{B}
\rput(-35,-27){$2$}
\cnode*(30,-30){2}{C}
\rput(35,-27){$3$}
\cnode(-50,-50){2}{D}
\rput(-50,-57){$Y_1$}
\cnode(-10,-50){2}{E}
\rput(5,-47){$4$}
\cnode*(10,-50){2}{F}
\rput(-10,-57){$Y_2$}
\cnode(50,-50){2}{G}
\rput(50,-57){$Y_5$}
\cnode(-10,-70){2}{H}
\rput(-10,-77){$Y_3$}
\cnode(30,-70){2}{K}
\rput(30,-77){$Y_4$}

\ncline{-}{A}{B}
\ncline{-}{A}{C}
\ncline{-}{A}{C}
\ncline{-}{B}{D}
\ncline{-}{B}{E}
\ncline{-}{C}{F}
\ncline{-}{C}{G}
\ncline{-}{F}{H}
\ncline{-}{F}{K}

\endpspicture
\caption{A tree in $\T_4$ of weight 
$(Y_1+Y_2)\,Y_1\,(Y_3+Y_4)\,Y_3$\,.}
\label{fig:2}
\end{figure}

Figure~\ref{fig:2} shows an example of a tree $T\in \T_4$ 
of weight $\wt(T)= 
(Y_1+Y_2)\,Y_1\,(Y_3+Y_4)\,Y_3$.
The vertices of $T$ are shown by 
black circles, and the added leaves of $\tilde T$ are shown by white circles.
The above recurrence relation for $P_r$ implies the following result.

\begin{proposition}
The parking polynomial $P_r$
equals the sum 
$$
P_r(Y_1,\dots,Y_r) = \sum_{T\in \T_r } \wt (T).
$$
\end{proposition}


\end{document}